\documentclass[twoside]{article}

\usepackage[usenames,dvipsnames]{color}
\usepackage{graphicx}

\usepackage{amssymb}

\usepackage{amsmath}
\usepackage{eufrak}
    \usepackage{ulem}
    \usepackage{chngcntr}
\usepackage{mathtools}
\usepackage{bm}

\usepackage{graphicx}
\usepackage{camnum}

\usepackage{epstopdf}

\newtheorem{remark}{Remark}
\counterwithin{table}{section}

\usepackage{tikz,pgfplots}
        \pgfplotsset{compat = 1.3}
        \pgfplotsset{minor grid style={dotted}} \pgfplotsset{major grid
        style={dashed}}

        \pgfplotsset{every x tick label/.append style={font=\footnotesize,
        yshift=0.25ex}}

        \pgfplotsset{every y tick label/.append
        style={font=\footnotesize, xshift=0.25ex}}

\usetikzlibrary{external}
        \definecolor{colorclassyorange}{rgb}{0.95000,0.32500,0.09800}
        \definecolor{colorchromeyellow}{rgb}{1.00000,0.6549,0}%
        \definecolor{colorpaleyellow}{rgb}{1.00000,0.8549,0.1}%

        \definecolor{colorclassyblue}{rgb}{0.00000,0.44706,0.74118}%
        \definecolor{colorpurple}{rgb}{0.49400,0.18400,0.55600}%
        \definecolor{colorfuschia}{rgb}{0.95039,0.0,0.95039}%
        \definecolor{colorlemongreen}{rgb}{0.6,0.8,0}%
        \definecolor{colorreal}{rgb}{0.92941,0.79412,0.12549}%
        \definecolor{colorimag}{rgb}{0.00000,0.49804,0.00000}%
        \definecolor{colorabs}{rgb}{1.00000,0.00000,0.00000}%

\newcommand{\ee}{{\mathrm e}}
\newcommand{\ii}{{\mathrm i}}

\newcommand{\FFF}{\mathcal{F}}
\newcommand{\AAA}{\mathcal{A}}
\newcommand{\DDD}{\mathcal{D}}
\newcommand{\WWW}{\mathcal{W}}

\newcommand{\LLL}{\mathcal{L}}

\newcommand{\dx}{\partial_x}
\newcommand{\dy}{\partial_y}

\newcommand{\MU}{\MM{\mu}}

\def\Frac#1#2{{\textstyle \frac{#1}{#2}}}

\font\Bbb=msbm10 scaled 1000
\font\sBbb=msbm7 scaled 1000
\font\gothic=eufm10 scaled 1100
\font\smallgothic=eufm7 scaled 1100
\font\curly=rsfs10 scaled 1000
\font\smallcurly=rsfs7 scaled 1000

\def\BB#1{\mbox{\Bbb #1}}              
\def\CC#1{\mathrm{#1}}                 
\def\GG#1{\mbox{\gothic #1}}           

\def\MM#1{\mbox{\boldmath$#1$\unboldmath}}

\newcommand{\Int}[4]{\int_{#2}^{#3}\!#4\,\mathrm{d}#1}

\newcommand{\ve}{\varepsilon}
\newcommand{\schr}{Schr\"odinger }

\newcommand{\TTT}{\mathcal{T}}

\newcommand{\norm}[2]{\left\| #2 \right\|_{#1}}

\def\OO#1{{\cal O}\left(#1\right)}

\numberwithin{equation}{section}

\usepackage[colorlinks,citecolor=red!50!black,linkcolor=red!50!black,urlcolor=red!50!black]{hyperref}
\usepackage{cleveref}
\begin{document}

\title{Sixth-order schemes for laser--matter interaction in the \schr equation}

\author{Pranav Singh\footnote{Mathematical Institute, Andrew Wiles Building, University of Oxford,
Radcliffe Observatory Quarter, Woodstock Rd, Oxford OX2 6GG, UK}} \maketitle

\vspace*{6pt} \noindent \textbf{AMS Mathematics Subject Classification:}
Primary 65M70, Secondary 35Q41, 65L05, 65F60

\vspace*{6pt} \noindent   \textbf{Keywords:} Schr\"odinger equation, laser potential, high order methods,
compact methods, splitting methods, time dependent potentials, Magnus
expansion

\begin{abstract}
Control of quantum systems via lasers has numerous applications that
require fast and accurate numerical solution of the \schr equation. In this
paper we present three strategies for extending any sixth-order scheme for
\schr equation with time-independent potential to a sixth-order method for
\schr equation with laser potential. As demonstrated via numerical
examples, these schemes prove effective in the atomic regime as well as the
semiclassical regime, and are a particularly appealing alternative to
time-ordered exponential splittings when the laser potential is highly
oscillatory or known only at specific points in time (on an equispaced
grid, for instance).

These schemes are derived by exploiting the linear in space form of the
time dependent potential under the dipole approximation (whereby
commutators in the Magnus expansion reduce to a simpler form), separating
the time step of numerical propagation from the issue of adequate
time-resolution of the laser field by keeping integrals intact in the
Magnus expansion, and eliminating terms with unfavourable structure via
carefully designed splittings.
\end{abstract}

\maketitle
\setcounter{section}{0}

\section{Introduction}
\label{sec:intro}
                In this paper we present a class of sixth-order numerical schemes for
                laser-matter interaction in the \schr equation under the dipole
                approximation,
                \begin{equation}\label{eq:schr}
                \ii \ve \partial_t u(\MM{x},t)=\big [-\ve^2 \Delta + V_0(\MM{x}) + \MM{e}(t)^\top  \MM{x}\big] u(\MM{x},t),\ u(\MM{x},0)=u_0(\MM{x}),
                \end{equation}
                where $\ t \geq 0$, $\MM{x}=(x_1,\ldots,x_n)\in\BB{R}^n$ and
                the laser term $\MM{e}(t)=(e_1(t),\ldots,e_n(t))$ is an $\BB{R}^n$ valued
                function of $t$. In accordance with convention, we will label the kinetic operator, $-\ve^2 \Delta$, its various scalings and certain closely related operators as $T$, while the potential operators such as $V_0(\MM{x}) + \MM{e}(t)^\top  \MM{x}$ and its various scalings will be labeled $W$.

                The parameter $\ve$ in \cref{eq:schr} acts like Planck's constant. In typical applications this
                parameter is $1$  when working in the atomic units and is very small, $0
                <\ve \ll 1$, when working in the semiclassical regime. The case $0<\ve \ll 1$ often appears through
                a rescaling of atomic units in applications involving heavier particles. To clarify the distinction, we refer to $\ve=1$ as the atomic
                {\em regime} and $0<\ve \ll 1$ as the semiclassical regime.
                 The methods developed
                in this paper are equally effective in both regimes.

                The direction of the
                time-dependent vector $\MM{e}(t)$ may be fixed in the case of linearly
                polarised light,
                \[ \MM{e}(t) = e(t) \hat{\MM{\mu}},\quad e(t)\in \BB{R},\ \hat{\MM{\mu}} \in \BB{R}^n,\ \norm{2}{\hat{\MM{\mu}}} =1, \]
                but we will treat it as a general time-dependent vector, which also covers the case of
                circular polarisation.

                \Cref{eq:schr} is a highly specialised case of the \schr equation with a
                time-dependent potential. Since lasers are among the most effective tools for
                controlling processes at the quantum scale  \cite{shapiro03qc} and  the dipole
                approximation is valid for a large range of applications, it is also a very
                important case that is encountered frequently in practice.

                Moreover, it often arises in some very challenging applications that require
                highly accurate but low cost numerical schemes. Agueny {\em et al.}  \cite{JanPetter2016}, for instance,
                consider a
                parameter sweep for a one-dimensional problem with spatial and temporal
                domains of size $10^5$ each (in atomic units). In such applications, the
                physical phenomenon under consideration is typically subtle and noticeable only
                over long temporal windows. This leads to a need for high order methods that are not only efficient enough
                to integrate over large temporal windows, but also
                accurate enough so that the error accumulated over large temporal windows stays reasonable.

                The solution of \cref{eq:schr} is also required in applications such the shaping of
                temporal profiles of lasers via optimal control where the numerical solutions
                for these equations are used repeatedly within an optimisation routine
                 \cite{AmstrupChirped,MeyerOptimal,CoudertOptimal}. Moreover, the laser field
                may be known only at specific times (such as on an equispaced grid)  and may be
                highly oscillatory in nature.

                A wide range of numerical schemes have  been designed for the case of
                time-dependent potentials
                 \cite{TalEzer1992,Peskin,SanzSerna96,TremblayCarrington04,KlaDimBrig2009,NdongEzerKosloff10,alvermann2011hocm,BIKS,BlanesCasasMurua17,blanes17quasimagnus,Schaefer,IKS18sinum,IKS18jcp}.
                Being general methods that are applicable to a broad class of time-dependent potentials, these methods can usually handle the case of \cref{eq:schr} well, which
                is often considered as an important example  \cite{TalEzer1992,Schaefer}. Nevertheless,
                these method often do not explicitly exploit the special structure of the
                laser potential under the dipole approximation.

                Among these, some of the notable classes are the family of Magnus--Lanczos methods \cite{KlaDimBrig2009,blanes17quasimagnus,IKS18sinum}
                that achieve high accuracies but require very small time step (or a large number of Lanczos iterations) and
                the class of Chebychev polynomial based methods  \cite{talezer84afs,NdongEzerKosloff10,Schaefer} that can allow the use of much larger time steps at the expense of unitarity.

                Yet another notable class is that of
                time-ordered exponential splittings  \cite{Suzuki93,Goldstein2004,Omelyan,SanzSerna96}, where a
                classical splitting for time-independent Hamiltonians is extended in a
                straightforward way for time-dependent Hamiltonians. However, since the
                time-knots where the potential is sampled are determined by the splitting
                coefficients, there is no flexibility in either (i) choosing specific knots
                if the potential is known only at specific times or (ii) using higher
                accuracy quadrature in the case of highly-oscillatory potentials.


                To overcome some of these limitations in the
                special case of \cref{eq:schr}, an efficient class of fourth-order numerical schemes
                were recently developed by resorting to exponential splittings of a
                fourth-order truncation of the Magnus expansion  \cite{Bandrauk2013,IKS18cpc}. These exploit the linearity
                (in space) of the time-dependent component of the potential and avoid discretisation of the integrals appearing in the Magnus expansion.
                While the sampling of potentials on an equispaced (temporal) grid
                can also be achieved effectively by following the approach of Schaeffer {\em et al.}  \cite{Schaefer}, the approach of Iserles {\em et al.}  \cite{IKS18cpc}
                was demonstrated to handle highly oscillatory potentials effectively very well.

                In this paper we extend the fourth-order techniques of Iserles {\em et al.}  \cite{IKS18cpc} to derive a class
                of sixth-order schemes that are specialised for the case of laser-matter
                interaction under the dipole approximation. These methods inherit many of the
                highly favourable properties of the fourth-order schemes.
                Namely,
                \begin{enumerate}
                \item the time integrals of $\MM{e}(t)$ are kept intact
                    till the very end, allowing their eventual
                    approximation via a variety of quadrature methods
                    depending on the nature of the laser pulse (such as
                    Newton--Cotes formulae for applications in optimal
                    control where the laser may be known at specific times,
                    or high-order Gauss--Legendre quadrature and Filon
                    quadrature for highly oscillatory lasers),
                \item the proposed schemes for laser potentials can be
                    implemented by extending existing  high-accuracy
                    implementations designed for \schr equation with
                    time-independent potentials at little to no extra
                    computational cost,
                \item and unlike Lanczos-based methods, which become
                    prohibitive for large time steps due to the large
                    spectral radius of the exponent, the cost of a single
                    time step of the scheme is either entirely independent
                    of the time step or grows mildly at worst, allowing the
                    use of large time steps. \end{enumerate}


                The simplification of the sixth-order Magnus expansion
                proves to be considerably more involved than the
                fourth-order case due to higher nested integrals and
                commutators. In particular, unlike the fourth-order Magnus
                expansion, it is not possible to reduce the sixth-order
                Magnus expansion to a commutator-free form. Moreover, this
                Magnus expansion features the gradient of the potential,
                which might be expensive or unavailable.

                The first option available to us for the exponentiation of
                this Magnus expansion is Lanczos iterations. Due to the
                structure of the simplified Magnus expansion, the cost of
                matrix-vector products in each Lanczos iteration is
                lower  \cite{IKS18sinum} compared to standard Magnus--Lanczos
                methods. The large number of iterations or small time steps
                typical of Magnus--Lanczos methods persist, however. A
                bigger advantage may be attained by splitting the
                exponential of the simplified Magnus expansion. Our
                simplified Magnus expansion, however,  is structurally more
                complex than the fourth-order expansion of Iserles {\em et
                al.}  \cite{IKS18cpc}, a hurdle also encountered by Goldstein
                \& Baye  \cite{Goldstein2004}, and requires entirely new
                splittings for effective exponentiation.

                For this purpose, we develop three different sixth-order
                exponential splittings. The first of these, and the closest
                to the approach of Iserles {\em et al.}  \cite{IKS18cpc},
                features a commutator and the gradient of the potential.
                The second is a specialised splitting that is free of
                commutators, while the third is free of commutators as well
                as the gradient of the potential.

                A common feature in these splittings is that the central
                exponent resembles the structure of the Hamiltonian and
                existing techniques for the \schr equation with
                time-independent potential can be readily applied for its
                exponentiation. The remaining exponentials turn out to be
                very inexpensive, thereby allowing us to extend any high
                accuracy method for time-independent potentials to the case
                of laser potentials at a very small additional cost.


                \subsection{Organization of the paper}
                We briefly revisit the fourth-order schemes of Iserles {\em
                et al.}  \cite{IKS18cpc} in \cref{sec:existing} before
                commencing the derivation of our sixth-order schemes. This
                proceeds in \cref{sec:magnus} by a simplification of the
                sixth-order Magnus expansion using commutator identities
                and integration-by-parts, eventually resulting in an
                expression that features a commutator and the gradient of
                the potential.

                The exponentiation of this expansion is addressed in
                \cref{sec:splitting}. Three exponential splitting
                strategies for this purpose are presented in
                \cref{sec:stranglanczos,sec:eliminatecomm,sec:gradfree}.
                These schemes fall under a common theme \cref{eq:common},
                where the central exponent needs to be approximated with an
                existing high-accuracy scheme for time-independent
                potentials. Three examples for this purpose -- a
                sixth-order classical splitting, a compact splitting and
                Lanczos approximation -- are described in \cref{sec:inner}.

                In \cref{sec:integrals} we discuss the approximation of
                integrals appearing in our schemes, while \cref{sec:expXYZ}
                describes the implementation of individual exponentials,
                completing the full description of the proposed schemes.
                Numerical examples for the proposed schemes are described
                in \cref{sec:numerics}, while our conclusions are
                summarised in \cref{sec:conclusions}.

\section{Existing fourth-order schemes}
\label{sec:existing} The \schr equation \cref{eq:schr} can be rewritten in
            the form
            \begin{equation}\label{eq:schr2}
            \partial_t u(\MM{x},t)= \AAA(\MM{x},t) u(\MM{x},t),\quad  \MM{x} \in  \BB{R}^n,\ t \geq 0,\ u(\MM{x},0)=u_0(\MM{x}),
            \end{equation}
            where $\AAA(\MM{x},t) = \ii \ve \Delta -\ii \ve^{-1} \left(
            V_0(\MM{x})+\MM{e}(t)^\top \MM{x}\right)$. In principle, a
            numerical solution of \cref{eq:schr2} can be given via the
            exponential of a truncated Magnus expansion, $\Theta_p$,
            \begin{equation}
            \label{eq:magapp}\MM{u}^{n+1} = \exp(\Theta_p(t_n+h,t_n))\, \MM{u}^n.\end{equation}
            In practice, approximating the exponential of the Magnus expansion can be quite challenging.
            Arguably the most popular approach for this purpose, the Lanczos iterations become very inefficient
            when moderate to long time steps are involved. This is because the superlinear accuracy of
            Lanczos approximation of the exponential is not achieved till the number of iterations  exceeds
            (roughly speaking) the spectral radius of the exponent  \cite{hochbruck97ksa}, effectively forcing the use of very small time steps
            (see \cref{sec:lanczos}).

            The fourth-order numerical schemes developed in Iserles {\em et
            al.}  \cite{IKS18cpc} overcome this difficulty by resorting to
            exponential splittings of a fourth-order truncation of the
            Magnus expansion,
            \[\Theta_2(t+h,t)=\Theta^{[1]}+\Theta^{[2]},\]
            where
            \begin{eqnarray*}
              \Theta^{[1]}&= &\Int{\zeta}{0}{h}{\AAA(t+\zeta)} =  \ii h \ve \Delta - \ii h \ve^{-1} \widetilde{V}(\MM{x},t,h),\\
              \Theta^{[2]}&= &-\Frac{1}{2}\Int{\zeta}{0}{h}{\Int{\xi}{0}{\zeta}{ [\AAA(t+\xi),\AAA(t+\zeta)]  }   }\\
              &=&- \Frac12 \left(\Int{\zeta}{0}{h}{\Int{\xi}{0}{\zeta}{
            \left[\MM{e}(t+\zeta) - \MM{e}(t+\xi) \right] } }\right)^\top [\Delta,
            \MM{x}] = -\MM{s}(t,h)^\top \nabla,
            \end{eqnarray*}
            since
            \begin{equation}
            \label{eq:DeltaX}
            [\Delta,\MM{x}]u = \sum_{j=1}^n \left( \partial_{x_j}^2 (\MM{x} u) - \MM{x} \partial_{x_j}^2 u \right) =  2 \nabla u,
            \end{equation}
            and
            \begin{align}
            \label{eq:Vtilde}
            \widetilde{V}(\MM{x},t,h) &= V_0(\MM{x}) + \MM{r}(t,h)^\top  \MM{x},\\
            \label{eq:r}\MM{r}(t,h) &= \Frac{1}{h} \Int{\zeta}{0}{h}{\MM{e}(t+\zeta)},\\
            \label{eq:s}\MM{s}(t,h) &= 2 \Int{\zeta}{0}{h}{
            \left(\zeta - \Frac{h}{2}\right)\MM{e}(t+\zeta)}.
            \end{align}

            The simplest of these is the scheme MaStBM
            (Magnus--Strang--Blanes--Moan),
            \begin{equation}
            \label{eq:MaStBM}
            \ee^{-\frac12 \MM{s}^\top  \nabla} \ee^{a_1 X} \ee^{b_1 Y}  \ee^{a_2 X} \ee^{b_2 Y} \ee^{a_3 X} \ee^{b_3 Y}  \ee^{a_4 X}  \ee^{b_3 Y}  \ee^{a_3 X} \ee^{b_2 Y} \ee^{a_2 X} \ee^{b_1 Y}  \ee^{a_1 X} \ee^{-\frac12 \MM{s}^\top  \nabla}, \end{equation}
            where $X = \ii h \ve \Delta$ and $Y = - \ii h \ve^{-1}
            \widetilde{V}(\MM{x},t,h)$. For
            the sake of brevity, we write $\MM{s}$ instead of
            $\MM{s}(t,h)$, suppressing $t$ and $h$.

            This approach splits the $\OO{h^3}$ term $\Theta^{[2]}$ from
            the Magnus expansion using Strang splitting and utilises a
            classical splitting of Blanes \& Moan  \cite{blanesandmoan} for
            the exponentiation of $\Theta^{[1]}$.

            Since $\MM{s}^\top \nabla$ commutes with the Laplacian, the outermost
            exponentials can be computed together in the form $\ee^{a_1 X -\frac12
            \MM{s}^\top  \nabla}$, without any additional cost compared to the classical
            splitting. The cost of computing a single exponential to arbitrary accuracy
            is independent of the time step since these are computed exactly via Fast
            Fourier Transforms (FFTs). A crucial advantage over time-ordered exponential
            splittings is that keeping the integrals intact in $\MM{r}$ and $\MM{s}$
            allows the sampling of the potential in more flexible ways.

\section{Simplification of the Magnus expansion}
\label{sec:magnus} In this section we present the first component in the
            derivation of our sixth-order schemes, which is the
            simplification of the sixth-order Magnus expansion,
            \[\Theta_4 = \Theta^{[1]}+\Theta^{[2]}+\Theta^{[3,1]}+\Theta^{[3,2]}+\Theta^{[4,1]}+\Theta^{[4,2]}+\Theta^{[4,3]},\]
             where the additional terms compared to $\Theta_2(t+h,t)$,
            \small
            \begin{eqnarray*}
              \Theta^{[3,1]}=& &\Frac{1}{12}\Int{\zeta}{0}{h}{\Int{\xi}{0}{\zeta}{ \Int{\chi}{0}{\zeta}{\left[\AAA(t+\xi),\left[ \AAA(t+\chi),\AAA(t+\zeta)\right]\right] } }}\\
              \Theta^{[3,2]}=& &\Frac{1}{4}\Int{\zeta}{0}{h}{\Int{\xi}{0}{\zeta}{ \Int{\chi}{0}{\xi}{\left[\left[\AAA(t+\chi), \AAA(t+\xi)\right],\AAA(t+\zeta)\right] } }},\\
                \Theta^{[4,1]}=&
                &-\Frac{1}{24}\Int{\zeta}{0}{h}{\Int{\xi}{0}{\zeta}{\Int{\chi}{0}{\zeta}{\Int{\nu}{0}{\chi}{\left[\AAA(t+\xi),
                  \left[\left[\AAA(t+\nu),\AAA(t+\chi)\right],\AAA(t+\zeta)\right]\right]}}}},\\
                \Theta^{[4,2]}=& &-\Frac{1}{24}\Int{\zeta}{0}{h}{\Int{\xi}{0}{\zeta}{\Int{\chi}{0}{\xi}{\Int{\nu}{0}{\xi}{\left[\left[\AAA(t+\chi),\left[\AAA(t+\nu),\AAA(t+\xi)\right]\right], \AAA(t+\zeta)\right] }}}}, \\
              \Theta^{[4,3]}=& &-\Frac{1}{8}\Int{\zeta}{0}{h}{\Int{\xi}{0}{\zeta}{\Int{\chi}{0}{\xi}{\Int{\nu}{0}{\chi}{\left[\left[\left[\AAA(t+\nu),
              \AAA(t+\chi)\right],\AAA(t+\xi)\right],\AAA(t+\zeta)\right]}}}},
              \end{eqnarray*}
            \normalsize are as specified by equation (4.18) in Iserles {\em
            et al.}  \cite{iserles00lgm}. We write $\Theta_4$, $\Theta^{[k]}$
            and $\Theta^{[k,j]}$ as a shorthand for $\Theta_4(t+h,t)$,
            $\Theta^{[k]}(t+h,t)$ and $\Theta^{[k,j]}(t+h,t)$,
            respectively, suppressing the dependence on $t$ and $h$ for
            brevity. Note that this is a {\em power-truncated Magnus
            expansion} where $\OO{h^7}$ terms have been discarded.

            \subsection{Simplification tools}
            In the simplification of commutators appearing in the
            sixth-order Magnus expansion, we will need the following
            commutator identities
            \begin{equation}
            \label{eq:ids} [\Delta, \MM{a}^\top \nabla] = 0, \quad  [\Delta, \MM{a}^\top \MM{x}] = 2 \MM{a}^\top \nabla, \quad [\MM{a}^\top \nabla, f] = \MM{a}^\top (\nabla f),
            \end{equation}
             where $\MM{a} \in \BB{C}^n$ and $f \in \CC{C_1}(\BB{R}^n;
            \BB{C})$. Further, for ease of computation, we write $\AAA(t)$
            in the form $\TTT + \WWW+ \LLL(t)$, where $\TTT = \ii \ve
            \Delta$, $\WWW= -\ii \ve^{-1} V_0(\MM{x})$ and $\LLL(t) = -\ii
            \ve^{-1} \MM{e}(t)^\top \MM{x}$ are linear differential
            operators that satisfy
            \[ [\TTT,\LLL(t)] = 2\MM{e}(t)^\top \nabla, \qquad [\WWW,\LLL(t)] = 0.\]

            The nested integrals in the Magnus expansion will, with one
            exception, be reduced to integrals over the interval,
            possessing a common form
%
%
            \begin{equation}
            \label{eq:mubernoulli} \MU_n(t,h) = \Int{\zeta}{0}{h}{\widetilde{B}_n(h, \zeta) \MM{e}(t+\zeta)},
            \end{equation}
             where $\widetilde{B}_n(h,\zeta) = h^n B_n(\zeta/h)$ is the
            n-th rescaled Bernoulli polynomial,
            \begin{equation}
            \label{eq:bernoulli}
            \widetilde{B}_0 = 1,\quad \widetilde{B}_1 = \zeta - \Frac12 h,\quad  \widetilde{B}_2 = \zeta^2 - h \zeta + \Frac16h^2,\quad  \widetilde{B}_3 = \zeta^3 - \Frac32 h \zeta^2 + \Frac12 h^2 \zeta.
            \end{equation}
            In this new notation,
            \[\MM{r}(t,h) = \Frac{1}{h}\MU_0(t,h), \qquad \MM{s}(t,h) = 2 \MU_1(t,h),\]
             which are similar to certain forms encountered in other
            Magnus-based approaches  \cite{IKS18sinum,IKS18jcp}.
            $\MM{r}(t,h),\MM{s}(t,h)$ and $\MU_n(t,h)$ will sometimes be
            abbreviated further to $\MM{r},\MM{s}$ and $\MU_n$,
            respectively, suppressing the dependence on $t$ and $h$.

            In the simplification of the nested integrals to integrals in
            the above form, we will frequently use a few identities derived
            via integration-by-parts.
            \begin{align}
            \label{eq:ibp1}\Int{\xi}{0}{\zeta}{\xi^n \Int{\chi}{0}{\xi}{\MM{e}(s+\chi)}}
            &=  \Int{\xi}{0}{\zeta}{\left(\zeta^{n+1}-\xi^{n+1}\right)
            \MM{e}(s+\xi)},\\
            \label{eq:ibp4} \Int{\xi}{0}{\zeta}{\Int{\chi}{0}{\xi}{\chi^n \MM{e}(s+\chi)}} &=  \Int{\xi}{0}{\zeta}{\left(\zeta \xi^n-\xi^{n+1}\right)
            \MM{e}(s+\xi)},\\
            \label{eq:ibp2}\Int{\xi}{0}{\zeta}{\xi^n \Int{\chi}{0}{\xi}{(\MM{e}(s+\xi)-\MM{e}(s+\chi))}} &= \Frac{1}{n+1}\Int{\xi}{0}{\zeta}{\left((n+2) \xi^{n+1} - \zeta^{n+1}\right)
            \MM{e}(s+\xi)},\\
            \label{eq:ibp3}
            \Int{\zeta}{0}{h}{\zeta \Int{\xi}{0}{\zeta}{ \xi\MM{e}(t+\xi)}} &= \Frac12\Int{\zeta}{0}{h}{\left(h^2 \zeta -\zeta^3\right)}.
            \end{align}

            \subsection{Simplification of $\Theta^{[3,1]}$}
            \label{sec:Th31} Using the above tools, the inner-most
            commutator in $\Theta^{[3,1]}$ is simplified as
            \begin{align}
            \nonumber  \left[ \AAA(t+\chi),\AAA(t+\zeta)\right] &= [\TTT + \WWW+ \LLL(t+\chi),\TTT + \WWW+ \LLL(t+\zeta)] \\
            \label{eq:grade2}& =  [\TTT,\LLL(t+\zeta)-\LLL(t+\chi)] = 2(\MM{e}(t+\zeta)-\MM{e}(t+\chi))^\top \nabla.
            \end{align}
             Using \cref{eq:grade2}, along with the following identities
            that result directly from \cref{eq:ids},
            \[ [\TTT, \MM{a}^\top \nabla] =0, \qquad [\WWW, \MM{a}^\top \nabla] = \ii \ve^{-1} \MM{a}^\top (\nabla V_0), \qquad [\LLL(t), \MM{a}^\top \nabla] = \ii \ve^{-1} \MM{a}^\top \MM{e}(t), \]
             we can simplify the full commutator in $\Theta^{[3,1]}$ to the
            form
            \begin{equation}
            \label{eq:grade3}
            \left[\AAA(t+\xi),\left[ \AAA(t+\chi),\AAA(t+\zeta)\right]\right]= 2 \ii \ve^{-1} (\MM{e}(t+\zeta)-\MM{e}(t+\chi))^\top \left(\nabla V_0+\MM{e}(t+\xi)\right).
            \end{equation}
             Combining this observation with \cref{eq:ibp2} under $n=1$, we
            can simplify
            \begin{eqnarray*}
            \Theta^{[3,1]} &= & \Frac{1}{6} \ii \ve^{-1}  \Int{\zeta}{0}{h}{\Int{\xi}{0}{\zeta}{ \Int{\chi}{0}{\zeta}{ (\MM{e}(t+\zeta)-\MM{e}(t+\chi))^\top \left((\nabla V_0)+\MM{e}(t+\xi)\right) }}}\\
            &=& \Frac{1}{6} \ii \ve^{-1} \left( \Int{\zeta}{0}{h}{\Int{\chi}{0}{\zeta}{ \zeta (\MM{e}(t+\zeta)-\MM{e}(t+\chi))}}\right)^\top (\nabla V_0)+c_{3,1}(t,h)\\
            &=& \Frac{1}{12} \ii \ve^{-1} \left(\Int{\zeta}{0}{h}{\left(3 \zeta^2 - h^2\right) \MM{e}(t+\zeta)}\right)^\top (\nabla V_0) +
            c_{3,1}(t,h),
            \end{eqnarray*}
            where
             $c_{3,1}(t,h)$ is a scalar whose simplification
            is confined to \cref{app:scalars}. In principle, this term can be ignored
            since it only results in a constant phase shift. Nevertheless, we carry it
            along for the sake of completeness.

            \subsection{Simplification of $\Theta^{[3,2]}$}
            The commutator in $\Theta^{[3,2]}$ is obtained from \cref{eq:grade3} by exchanging $\chi$ and $\zeta$,
            \begin{eqnarray*}
            \left[\left[\AAA(t+\chi), \AAA(t+\xi)\right],\AAA(t+\zeta)\right] &=& \left[\AAA(t+\zeta),\left[ \AAA(t+\xi),\AAA(t+\chi)\right]\right] \\
            &=& 2 \ii \ve^{-1} (\MM{e}(t+\chi)-\MM{e}(t+\xi))^\top \left(\nabla
            V_0+\MM{e}(t+\zeta)\right).
            \end{eqnarray*}
            The simplification of $\Theta^{[3,2]}$ results by using \cref{eq:ibp2} once
            under $n=0$ (for the two inner integrals), followed by an application of
            \cref{eq:ibp1,eq:ibp4} under $n=1$,
            \begin{eqnarray*}
            \Theta^{[3,2]} &= &\Frac{1}{2}\ii \ve^{-1} \left(\Int{\zeta}{0}{h}{\Int{\xi}{0}{\zeta}{ \Int{\chi}{0}{\xi}{ (\MM{e}(t+\chi)-\MM{e}(t+\xi))}}}\right)^\top (\nabla V_0) + c_{3,2}(t,h)\\
            &=& -\ii \ve^{-1}  \left(\Int{\zeta}{0}{h}{\Int{\xi}{0}{\zeta}{\left(\xi - \Frac{\zeta}{2}\right)\MM{e}(t+\xi)} }\right)^\top (\nabla V_0) + c_{3,2}(t,h)\\
            &=&\Frac14 \ii \ve^{-1} \left(\Int{\zeta}{0}{h}{\left(3\zeta^2 + h^2 - 4 h \zeta\right) \MM{e}(t+\zeta)}\right)^\top
            (\nabla V_0) + c_{3,2}(t,h).
            \end{eqnarray*}
            The simplification of the scalar $c_{3,2}(t,h)$ is, once again, confined to
            \cref{app:scalars}.

            Putting these results together, we find
            \begin{equation}
            \label{eq:Th3}
            \Theta^{[3]} :=\Theta^{[3,1]}+\Theta^{[3,2]} = \ii \ve^{-1} \MM{q}(t,h)^\top (\nabla V_0) +c(t,h),
            \end{equation}
             where, using the $\MU_n$ notation of \cref{eq:mubernoulli},
            \begin{equation}
            \label{eq:q}
            \MM{q}(t,h) = \Int{\zeta}{0}{h}{\left(\zeta^2 - h \zeta + \Frac16h^2\right) \MM{e}(t+\zeta)}  = \MU_2(t,h),
            \end{equation}
             and $c(t,h)= c_{3,1}(t,h)+c_{3,2}(t,h)$ is
            \begin{align}
            \nonumber
            c(t,h) &= \ii \ve^{-1}\left( 2 \Int{\zeta}{0}{h}{\zeta \MM{e}(t+\zeta)^\top \Int{\xi}{0}{\zeta}{\MM{e}(t+\xi)}} - \left(\Int{\zeta}{0}{h}{\MM{e}(t+\zeta)}\right)^\top
            \left(\Int{\zeta}{0}{h}{\zeta \MM{e}(t+\zeta)}\right) \right. \\
            \label{eq:c}& \qquad \qquad\left.  - \Frac{1}{6} h \left(\Int{\zeta}{0}{h}{\MM{e}(t+\zeta)}\right)^2\right),
            \end{align}
             from \cref{app:scalars}.

            \subsection{Simplification of $\Theta^{[4,1]}$}
            We simplify the commutator in $\Theta^{[4,1]}$ starting from the result of
            \cref{eq:grade3},
            \begin{eqnarray*}
             & &\left[\AAA(t+\xi),
                  \left[\left[\AAA(t+\nu),\AAA(t+\chi)\right],\AAA(t+\zeta)\right]\right]  \\
            & & \qquad = - 2 \ii \ve^{-1} \left[ \TTT+\WWW+\LLL(t+\xi), (\MM{e}(t+\chi)-\MM{e}(t+\nu))^\top \left(\nabla V_0+\MM{e}(t+\zeta)\right)\right]\\
            & & \qquad = 2  \left[\Delta, (\MM{e}(t+\chi)-\MM{e}(t+\nu))^\top (\nabla V_0)\right].\\
            \end{eqnarray*}
            We avoid simplifying the above commutator further since doing
            so does not give us any computational advantage.  The component
            simplifies to
            \begin{equation}
            \label{eq:Th41}
            \Theta^{[4,1]} =  \left[\Delta, \MM{p}_1(t,h)^\top (\nabla V_0)\right],
            \end{equation}
             where the four times nested integral is reduced to an integral
            over an interval,
            \begin{align}
            \nonumber \MM{p}_1(t,h)&= -\Frac{1}{12}\Int{\zeta}{0}{h}{\Int{\xi}{0}{\zeta}{\Int{\chi}{0}{\zeta}{\Int{\nu}{0}{\chi}{  (\MM{e}(t+\chi)-\MM{e}(t+\nu))}}}}\\
            \nonumber&= -\Frac{1}{6}\Int{\zeta}{0}{h}{\zeta \Int{\chi}{0}{\zeta}{(\chi-\Frac{\zeta}{2}) \MM{e}(t+\chi)}}\\
            \nonumber &= -\Frac{1}{12}\Int{\zeta}{0}{h}{\left(h^2 \zeta-\zeta^3 \zeta\right) \MM{e}(t+\zeta)}  +\Frac{1}{36}\Int{\zeta}{0}{h}{\left(h^{3}-\zeta^{3}\right) \MM{e}(t+\zeta)} \\
            \label{eq:p1}&= \Frac{1}{36} \Int{\zeta}{0}{h}{\left(2 \zeta^3 - 3 h^2 \zeta+ h^{3}\right) \MM{e}(t+\zeta)}.
            \end{align}
             using \cref{eq:ibp2} under $n=0$, \cref{eq:ibp3} and
            \cref{eq:ibp1} under $n=2$.

            \subsection{Simplification of $\Theta^{[4,2]}$}
            Reusing the workings of $\Theta^{[4,1]}$, the commutator in
            $\Theta^{[4,2]}$ is
            \[\left[\left[\AAA(t+\chi),\left[\AAA(t+\nu),\AAA(t+\xi)\right]\right], \AAA(t+\zeta)\right] = 2  \left[\Delta, (\MM{e}(t+\xi)-\MM{e}(t+\nu))^\top (\nabla V_0)\right].\]
             Integrating the occurrence of this commutator in
            $\Theta^{[4,2]}$, we find
            \begin{equation}
            \label{eq:Th42}
            \Theta^{[4,2]} =  \left[\Delta, \MM{p}_2(t,h)^\top (\nabla V_0)\right],
            \end{equation}
             and the component simplifies to
            \begin{align}
            \nonumber \MM{p}_2(t,h)&= -\Frac{1}{12}\Int{\zeta}{0}{h}{\Int{\xi}{0}{\zeta}{\Int{\chi}{0}{\xi}{ \Int{\nu}{0}{\xi}{(\MM{e}(t+\xi)-\MM{e}(t+\nu))}}}}\\
            \nonumber &= -\Frac{1}{12}\Int{\zeta}{0}{h}{\Int{\xi}{0}{\zeta}{\Int{\nu}{0}{\xi}{ \xi (\MM{e}(t+\xi)-\MM{e}(t+\nu))}}}\\
            \nonumber &= -\Frac{1}{8}\Int{\zeta}{0}{h}{\Int{\xi}{0}{\zeta}{\xi^2\MM{e}(t+\xi)}}+\Frac{1}{24}\Int{\zeta}{0}{h}{\zeta^2 \Int{\xi}{0}{\zeta}{\MM{e}(t+\xi)}}\\
            \nonumber &= - \Frac{1}{8}\Int{\zeta}{0}{h}{\left(h \zeta^2 - \zeta^3\right) \MM{e}(t+\zeta)} +\Frac{1}{72} \Int{\zeta}{0}{h}{\left(h^3-\zeta^3\right) \MM{e}(t+\zeta)}\\
            \label{eq:p2}& = \Frac{1}{72}\Int{\zeta}{0}{h}{\left(8 \zeta^3 -9 h \zeta^2 + h^3\right) \MM{e}(t+\zeta)}
            \end{align}
             by using \cref{eq:ibp2} under $n=1$, and
            \cref{eq:ibp1,eq:ibp4} under $n=2$.

            \subsection{Simplification of $\Theta^{[4,3]}$}
            The last remaining term is $\Theta^{[4,3]}$, which features the
            commutator
            \[
            \left[\left[\left[\AAA(t+\nu),
              \AAA(t+\chi)\right],\AAA(t+\xi)\right],\AAA(t+\zeta)\right] = -2 \left[\Delta,  (\MM{e}(t+\chi)-\MM{e}(t+\nu))^\top (\nabla
            V_0)\right],\] that simplifies to
            \begin{equation}
            \label{eq:Th43}
            \Theta^{[4,3]} =  \left[\Delta, \MM{p}_3(t,h)^\top (\nabla V_0)\right],
            \end{equation}
             where
            \begin{align}
            \nonumber \MM{p}_3(t,h)&= \Frac{1}{4}\Int{\zeta}{0}{h}{\Int{\xi}{0}{\zeta}{\Int{\chi}{0}{\xi}{\Int{\nu}{0}{\chi}{ (\MM{e}(t+\chi)-\MM{e}(t+\nu))}}}}\\
            \nonumber&= \Frac{1}{4}\Int{\zeta}{0}{h}{\Int{\xi}{0}{\zeta}{\Int{\chi}{0}{\xi}{(2\chi - \xi)\MM{e}(t+\chi)}}}\\
            \nonumber&= \Frac{1}{2}\Int{\zeta}{0}{h}{\Int{\xi}{0}{\zeta}{(\zeta \xi - \xi^2)e(t+\xi)}}
            -\Frac{1}{8}\Int{\zeta}{0}{h}{\Int{\xi}{0}{\zeta}{(\zeta^2 - \xi^2)e(t+\xi)}}
            \\
            \nonumber&=  \Int{\zeta}{0}{h}{\left[\Frac14\left(h^2 \zeta - \zeta^3\right) -\Frac38 \left(h \zeta^2 - \zeta^3\right) - \Frac{1}{24} \left(h^3-\zeta^3\right) \right]\MM{e}(t+\zeta)}
            \\
            \label{eq:p3}&= \Frac{1}{24}\Int{\zeta}{0}{h}{\left(4 \zeta^3 - 9 h \zeta^2 + 6 h^2 \zeta  - h^3\right) \MM{e}(t+\zeta)},
            \end{align}
             by using \cref{eq:ibp2} under $n=0$, \cref{eq:ibp1,eq:ibp4}
            under $n=1$, \cref{eq:ibp3} and, finally,
            \cref{eq:ibp1,eq:ibp4} under $n=2$.

            Putting these together,
            \begin{equation}
            \label{eq:Th4}
            \Theta^{[4]} :=\Theta^{[4,1]}+\Theta^{[4,2]}+\Theta^{[4,3]} = \left[\Delta, \MM{p}(t,h)^\top (\nabla V_0)\right],
            \end{equation}
             where, using \cref{eq:p1,eq:p2,eq:p3},
            \begin{equation}
            \label{eq:p}
            \MM{p}(t,h) = \Frac{1}{3} \Int{\zeta}{0}{h}{\left(\zeta^3 - \Frac32 h \zeta^2 - \Frac12 h^3\right) \MM{e}(t+\zeta)} = \Frac13 \MU_3(t,h).
            \end{equation}

            \subsection{The sixth-order Magnus expansion}
            \label{sec:Mag} The derivation of the simplified sixth-order
            Magnus expansion is completed by adding \cref{eq:Th3,eq:Th4} to
            the simplified fourth-order Magnus expansion  \cite{IKS18cpc}
            described in \cref{sec:existing}. Collecting these results, the
            sixth-order Magnus expansion is
            \begin{equation}
            \label{eq:Mag} \Theta_4 = \ii h \ve \Delta - \ii h \ve^{-1} \widetilde{V}  - \MM{s}^\top \nabla + \ii \ve^{-1} \MM{q}^\top (\nabla V_0)  +  \left[\Delta, \MM{p}^\top (\nabla V_0)\right] +c,
            \end{equation}
             where
            \[\tag{\ref{eq:Vtilde}}
            \widetilde{V}(\MM{x};t,h) = V_0(\MM{x}) + \MM{r}(t,h)^\top
            \MM{x},\] the coefficients $\MM{r},\MM{s},\MM{q},\MM{p}\in
            \BB{R}^n$ are given by
            \[\tag{\ref{eq:r},\ref{eq:s},\ref{eq:q},\ref{eq:p}} \MM{r} = \Frac{1}{h}\MU_0, \quad \MM{s} = 2 \MU_1,
            \quad  \MM{q}  = \MU_2, \quad \MM{p}  = \Frac{1}{3} \MU_3,\]
            and $c \in \BB{C}$ is given by \cref{eq:c}.
            %
            For the purpose of completion, we recall the definitions
            \begin{equation}
            \tag{\ref{eq:mubernoulli}} \MU_n(t,h) = \Int{\zeta}{0}{h}{\widetilde{B}_n(h, \zeta) \MM{e}(t+\zeta)},
            \end{equation}
             and
            \[\tag{\ref{eq:bernoulli}}
            \widetilde{B}_0 = 1,\quad \widetilde{B}_1 = \zeta - \Frac12
            h,\quad \widetilde{B}_2 = \zeta^2 - h \zeta + \Frac16h^2,\quad
            \widetilde{B}_3 = \zeta^3 - \Frac32 h \zeta^2 + \Frac12 h^2
            \zeta .\] Traditional Magnus--Lanczos methods directly
            approximate $\exp(\Theta)\MM{v}$ via Lanczos iterations, which
            involve computation of matrix-vector products of the form
            $\Theta \MM{v}$. These matrix-vector products are costly to
            compute since they involve nested commutators. Many approaches
            for reducing the nested commutators have been developed over
            the years  \cite{mko99fla,blanes00iho}.

            Approximating the exponential of \cref{eq:Mag}, $\exp(\Theta_4) \MM{v}$, via Lanczos iterations
            is an appealing option that should prove cheaper  \cite{IKS18sinum} than
            existing Magnus--Lanczos schemes since \cref{eq:Mag} involves only
            one commutator, leading to a lower cost of matrix-vector products. Moreoever, \cref{eq:Mag} preserves integrals which leads to high accuracy for oscillatory (in time) potentials  \cite{IKS18sinum}. These advantages must naturally be weighed against the
            requirement for $\nabla V_0$.
            Since Lanczos approximation of the exponential is relatively straightforward and well understood, however, this approach will not be elaborated in detail.
            Certain aspects of Lanczos iterations are touched upon in \cref{sec:lanczos} briefly.

            An appealing alternative to Lanczos approximation is to split
            the exponential of the Magnus expansion, effective strategies
            for which have been developed up to order four  \cite{IKS18cpc}.
            The sixth-order Magnus expansion, in the context of dipole
            approximation in laser-matter interaction, is structurally
            different from the fourth-order Magnus expansion, however, due
            to the appearance of the commutator $\left[\Delta, \MM{p}^\top
            (\nabla V_0)\right]$, which presents many additional
            difficulties in the exponentiation. Moreover, unlike the
            fourth-order Magnus expansion, $\Theta_4$ features the gradient
            of the potential, $\nabla V_0$, which might be unavailable.
            Specialised exponential splittings are required, therefore, if
            we wish to do without either the commutator or the gradient of
            the potential. The development of such schemes is pursued in
            \cref{sec:splitting}.

            \subsection{Sizes of components}
            \label{sec:evenodd} An essential ingredient in an effective exponential
            splitting of the Magnus expansion is a good estimate of the sizes of various
            components of in terms of the time step, $h$.

            For instance, since $\widetilde{B}_1(h,\zeta) = \OO{h}$ for
            $\zeta = \OO{h}$, one might expect that $\MM{s} = \MU_1 =
            \Int{\zeta}{0}{h}{\widetilde{B}_1(h,\zeta)\MM{e}(t+\zeta)}
            $ scales as $\OO{h^2}$. However, this overestimates the size of
            $\MM{s}$ when $\MM{e}$ is analytic. To see this, consider any
            $f$ such that
            \[ \Int{\zeta}{0}{h}{f(h, \zeta)} = 0, \quad \text{and}\quad f(h,\zeta) = \OO{h^n}\ \text{for}\ \zeta = \OO{h}. \]
             Expanding $\MM{e}$ at $\hat{t} \in [t,t+h]$ using Taylor
            expansion,
            \begin{equation}
            \label{eq:gainofpower}
            \Int{\zeta}{0}{h}{f(h, \zeta) \MM{e}(t+\zeta)} =  \sum_{k=1}^\infty \MM{e}^{(k)}(\hat{t}) \Int{\zeta}{0}{h}{f(h,
            \zeta) (t-\hat{t}+\zeta)^k} =  \OO{h^{n+2}},\end{equation}
             we find that the
            $k=0$ term vanishes, since $\Int{\zeta}{0}{h}{f(h, \zeta)} = 0$.
            Consequently, the first non-vanishing term in the expansion is the $k=1$
            term, and $\Int{\zeta}{0}{h}{f(h, \zeta) \MM{e}(t+\zeta)}$ ends up being $\OO{h^{n+2}}$ in this special case, instead
            of $\OO{h^{n+1}}$ which would usually be expected when only working under the
            assumption $f=\OO{h^n}$.

            Since integrals of the (re-scaled) Bernoulli polynomials
            vanish,
            \[\Int{\zeta}{0}{h}{\widetilde{B}_n(h, \zeta)} = 0\quad \text{for}\ n\geq 1, \quad
            \text{and} \quad \widetilde{B}_n(h,\zeta) = \OO{h^n}\ \text{for}\ \zeta =
            \OO{h},\]  we conclude that
             \[\MU_n = \OO{h^{n+2}},\ n\geq 1.\]
            Consequently, we expect that
            \begin{equation*}
            \label{eq:powersofhsimple}
            \MM{r}(t,h) = \OO{h},\quad \MM{s}(t,h) = \OO{h^3}, \quad \MM{q}(t,h) = \OO{h^4}, \quad \MM{p}(t,h) = \OO{h^5}.
            \end{equation*}

            This estimate is in line with the standard analysis presented
            in for Magnus expansions  \cite{iserles99ots,iserles00lgm} in a
            much more general setting and is not surprising. The $\OO{h^4}$
            size of $\MM{q}$ estimated via this analysis, however, turns
            out to be too large for our purposes, causing many difficulties
            in the design and analysis of the exponential splittings. In
            the case of compact splittings described in \cref{sec:compact},
            for instance, this would normally force us to compute $\nabla
            (\MM{q}^\top (\nabla V_0))$, which involves computing mixed
            derivatives of the potential $V_0$.

            This is remedied easily by noting that
            \[\widetilde{B}_2(h,\zeta) = \zeta^2 - h \zeta - \Frac16 h^2 = \left(\zeta-\Frac{h}{2}\right)^2 - \Frac{5}{12} h^2\]
             is even around $\Frac{h}{2}$. Expanding $\MM{e}$ at the
            midpoint of the interval, $\hat{t} = t+\Frac{h}{2}$, we find
            that the $k=1$ term in the Taylor expansion
            \cref{eq:gainofpower} for $\MU_2$,
            \[ \Int{\zeta}{0}{h}{\widetilde{B}_2(h, \zeta) (t-\hat{t}+\zeta)} =
            \Int{\zeta}{0}{h}{(\zeta-\Frac{h}{2})^3} - \Frac{5}{12} h^2
            \Int{\zeta}{0}{h}{(\zeta-\Frac{h}{2})} = 0, \]  also vanishes
            in addition to the $k=0$ term. Consequently, $\MM{q} = \MU_2$
            is $\OO{h^5}$, not $\OO{h^4}$ as suggested by standard
            analysis.

            Summarising our observations,
            \begin{equation}
            \label{eq:powersofh}
            \MM{r}(t,h) = \OO{h},\quad \MM{s}(t,h) = \OO{h^3}, \quad \MM{q}(t,h) = \OO{h^5}, \quad \MM{p}(t,h) = \OO{h^5}.
            \end{equation}

\section{Exponential splittings for the Magnus expansion}
\label{sec:splitting}
            As mentioned previously  in \cref{sec:intro,sec:existing}, the numerical
            exponentiation of the Magnus expansion, \cref{eq:Mag}, is incredibly costly
            unless split in a clever fashion.

            The common theme among our splittings
            \begin{equation}\label{eq:common} \ee^{\Theta_4(t+h,t)} = \ee^{\frac12 L}\ee^{\frac12 C}\ee^{T+W}\ee^{\frac12 C} \ee^{\frac12 L} + \OO{h^7},\end{equation}
            will be that they express the exponential of the Magnus expansion up to order six accuracy in terms of
            products of five or three (under L = 0) exponentials. While the forms of $C$ and $L$ will vary in the different splittings,
            as will the exact expressions of $T$ and $W$, what remains common is that $T$
            is a modified kinetic term and $W$ is a modified potential term (for instance,  $T_1 = h \TTT -  \MM{s}^\top \nabla$ in \cref{eq:parts1}). In particular, the structure is chosen to ensure that the separate exponentials of $T$ and $W$ are very
            inexpensive to compute exactly. Consequently, the inner-most exponential $\ee^{T+W}$
            can be approximated very efficiently via existing exponential splitting schemes for \schr equations with time-independent potentials.

            In \cref{sec:stranglanczos,sec:eliminatecomm,sec:gradfree}, we will develop
            three different sixth-order exponential splittings,
            \cref{eq:stranglanczos,eq:commfree,eq:gradfree}, that prescribe to the common
            form \cref{eq:common}. To fully describe concrete examples of these schemes,
            we consider two types of sixth-order splittings for approximating $\ee^{T+W}$
            in \cref{sec:inner}

\subsection{Schemes featuring a commutator}
\label{sec:stranglanczos} Our first sixth-order splitting is obtained via a
            Strang splitting of the Magnus expansion, \cref{eq:Mag},
            \begin{equation}
            \label{eq:stranglanczos}
            \tag{S1}
            \ee^{\frac{1}{2} \left[\Delta, \MM{p}^\top (\nabla V_0)\right]}
            \ee^{(\ii h \ve \Delta   - \MM{s}^\top \nabla) + (- \ii h \ve^{-1}  \widetilde{V} + \ii \ve^{-1} \MM{q}^\top (\nabla V_0) +c)}
            \ee^{\frac{1}{2} \left[\Delta, \MM{p}^\top (\nabla V_0)\right]}\,,
            \end{equation}
             where, in the context of \cref{eq:common},
            \begin{align}
            \nonumber T_1 &= \ii h \ve \Delta   - \MM{s}^\top \nabla, & W_1 &= - \ii h \ve^{-1}  \widetilde{V} + \ii \ve^{-1} \MM{q}^\top (\nabla V_0) +c, \\
            \label{eq:parts1} C_1 &=  \left[\Delta, \MM{p}^\top (\nabla V_0)\right], & \ L_1 &= 0,
            \end{align}
             and where the subscript in $T_1,W_1,C_1$ and $L_1$ indicates
            that these describe the first of the three classes of
            splittings presented in this manuscript.

            \subsubsection{Sixth-order accuracy of the scheme}
            Writing $\Theta_4(t+h,t)$ as $T_1+W_1+C_1$, the Strang
            splitting,
            \begin{equation}\label{eq:strang} \ee^{\Theta_4(t+h,t)} = \ee^{\frac12 C_1} \ee^{T_1+W_1} \ee^{\frac12  C_1} +
            \OO{h^7},
            \end{equation}
             turns out to be an order six splitting. To see this, recall
            the symmetric Baker--Campbell--Hausdorff (sBCH) formula, \small
            \begin{equation}\label{eq:sBCH} \ee^{\frac12 A} \ee^{B} \ee^{\frac12 A} = \ee^{\mathrm{sBCH}(A,B)}, \ \mathrm{sBCH}(A,B) = A+ B - \left(\Frac{1}{24}[[B,A],A] +
            \Frac{1}{12} [[B,A],B] \right) + \mathrm{h.o.t.}\end{equation}
            \normalsize
            Recall from \cref{eq:powersofh}  that $\MM{s} = \OO{h^3}$, $\MM{q} = \OO{h^5}$ and $\MM{p} =
            \OO{h^5}$, so that $T_1$ and $W_1$ are $\OO{h}$,
            while $C_1$ is $\OO{h^5}$. Thus,
            $\ee^{\frac12 C_1} \ee^{T_1+W_1} \ee^{\frac12 C_1}$ differs from $\ee^{ T_1 + W_1
            +C_1}$ by $\Frac{1}{24}[[C_1,T_1 + W_1],T_1 + W_1]$ (or smaller terms), which happens to be $\OO{h^7}$.

\subsection{Eliminating commutators}
\label{sec:eliminatecomm} Although the commutator in the splitting
                \cref{eq:stranglanczos} tends to be fairly benign (see
                \cref{sec:Z}), it can potentially be problematic at very
                large time steps. In this section we develop a specialised
                splitting  that overcomes this limitation,
                \begin{equation}
                \label{eq:commfree}
                \tag{S2}
                \ee^{- 6  h^{-2} \MM{p}^\top \nabla}
                \ee^{(\ii h \ve \Delta -  \tilde{\MM{s}}^\top \nabla)+(- \ii h \ve^{-1}  \widetilde{V} + \ii \ve^{-1} \MM{q}^\top (\nabla V_0) +c)}
                \ee^{- 6  h^{-2} \MM{p}^\top \nabla},
                \end{equation}
                 where
                \begin{equation}
                \label{eq:stilde}
                \tilde{\MM{s}} = \MM{s} -12  h^{-2}\MM{p}.
                \end{equation}
                 In the context of \cref{eq:common},
                \begin{align}
                \nonumber T_2 &= T_1-C_2 = \ii h \ve \Delta   - \tilde{\MM{s}}^\top \nabla, & W_2 &= W_1 = - \ii h \ve^{-1}  \widetilde{V} + \ii \ve^{-1} \MM{q}^\top (\nabla V_0) +c, \\
                \label{eq:parts2} C_2 &= - 12  h^{-2} \MM{p}^\top \nabla, & \ L_2 &= 0,
                \end{align}
                 where $T_2$ is a slight perturbation of $T_1$ but
                maintains the same structure, while $W_2$ is identical to
                $W_1$  defined in \cref{eq:parts1}.

                \begin{remark}
                \label{rmk:combineCT} Crucially, $C_2$ is free of commutators and, in fact,
                commutes with $T_2$. Consequently, in an exponential splitting of
                $\exp(T_2+W_2)$ where the outermost exponent happens to be $\ee^{a_1 T_2}$,
                the exponential $\ee^{C_2}$ can be combined with it (see \cref{eq:common})
                so that we only need to compute $\ee^{a_1 T_2 + C_2}$. This is the case for
                both splittings of $\exp(T_2+W_2)$ that are presented in \cref{sec:inner}.
                When combined with such splittings, our second class of sixth-order
                splittings for laser potentials, \cref{eq:commfree}, features no additional
                exponential in comparison to existing sixth-order schemes for
                time-independent potentials.
                \end{remark}
                \subsubsection{Derivation of the scheme}
                We start by letting
                \begin{equation} \label{eq:U} C_2 = \lambda h^{-2} \MM{p}^\top \nabla
                \end{equation}
                for some $\lambda \in \BB{R}$ to be determined, and attempt
                to express the exponential of the Magnus expansion in the
                form \cref{eq:commfree},
                \begin{equation}
                \label{eq:commfreeattempt}
                \exp(\Theta_4) =  \ee^{\frac12 C_2} \ee^{(T_1-C_2)+W_1}  \ee^{\frac12 C_2} + \OO{h^7}.
                \end{equation}
                 In order to express the right side of
                \cref{eq:commfreeattempt} as a single exponential, we use
                the sBCH formula \cref{eq:sBCH} up to an accuracy of
                $\OO{h^7}$ with $A=C_2$ and $B=(T_1 - C_2) +W_1$, \small
                \begin{eqnarray*}
                & &\mathrm{sBCH}(C_2,(T_1-C_2)+W_1) \\
                & &\qquad = T_1+W_1 - \left(\Frac{1}{24}[[T_1-C_2+W_1,C_2],C_2] + \Frac{1}{12} [[T_1-C_2+W_1,C_2],T_1-C_2+W_1] \right) \\
                & &\qquad = T_1+W_1- \left(\Frac{1}{24}[[W_1,C_2],C_2] + \Frac{1}{12} [[W_1,C_2],T_1-C_2+W_1] \right)\\
                & &\qquad = T_1+W_1- \Frac{1}{12} [[W_1,C_2],T_1+W_1].
                \end{eqnarray*}
                \normalsize Note that, since $\MM{p}=\OO{h^5}$, $C_2$ scales as $\OO{h^3}$.
                Consequently, grade five commutators (which feature five occurrences of $A$
                and $B$) involving even a single occurrence of $C_2$ are $\OO{h^7}$ or
                smaller and can be ignored. Thus it suffices to truncate the sBCH at grade
                three. We have also utilised the fact that $T_1$ (and $C_2$ itself)
                commutes with $C_2$ and drops out of the inner commutators. The grade three
                commutator $[[W_1,C_2],C_2]$ is also $\OO{h^7}$ due to two occurrences of
                $C_2$, and can be ignored. In the only remaining non-trivial term,
                $[[W_1,C_2],T_1-C_2+W_1]$, the component $[[W_1,C_2],C_2]$ can once again
                be ignored due to size.

                At this stage we compute the inner commutator up to
                $\OO{h^7}$ accuracy,
                \[[W_1,C_2] = -\lambda h^{-2} [\MM{p}^\top \nabla,   - \ii h
                \ve^{-1} \widetilde{V} + \ii \ve^{-1} \MM{q}^\top (\nabla V_0) +c ] =  \ii
                \lambda h^{-1} \ve^{-1}  [\MM{p}^\top \nabla,
                 \widetilde{V}] + \OO{h^7}, \] where the term involving
                 both $\MM{p}$ and $\MM{q}$ , which are both $\OO{h^5}$, is
                too small and can be ignored. Thus, $[W_1,C_2] = \ii
                \lambda h^{-1} \ve^{-1} \MM{p}^\top (\nabla  \widetilde{V})
                + \OO{h^7}$ using \cref{eq:ids}. This is a function (or a
                multiplication operator) and, consequently, commutes with
                $W_1$. Thus, the only relevant term of
                $[[W_1,C_2],T_1-C_2+W_1]$ is $[[W_1,C_2],T_1]$, i.e.
                \[[[W_1,C_2],T_1-C_2+W_1] =  \lambda [\Delta, \MM{p}^\top (\nabla  \widetilde{V})] +\OO{h^7}. \]
                 Under the choice of $\lambda = -12$,
                \begin{equation*}
                \mathrm{sBCH}(C_2,(T_1-C_2)+W_1) = T_1 + W_1 + C_1 + \OO{h^7},
                \end{equation*}
                 since,
                \[ - \Frac{1}{12} \lambda  [\Delta, \MM{p}^\top (\nabla \widetilde{V})] = [\Delta, \MM{p}^\top (\nabla V_0)] +  [\Delta, \MM{p}^\top \MM{r}] =  [\Delta, \MM{p}^\top (\nabla V_0)] = C_1.\]
                 In other words, since $\Theta_4 = T_1 + W_1 + C_1 $,
                \[ \ee^{\frac12 C_2} \ee^{T_2+W_2}  \ee^{\frac12 C_2}  = \ee^{\Theta_4} + \OO{h^7},\]
                 is a sixth-order splitting for \cref{eq:Mag}.

\subsection{Eliminating gradients of the potential}
\label{sec:gradfree} In this section we develop a specialised splitting
            that requires neither commutators nor the gradient of the
            potential, $\nabla V_0$. This sixth-order splitting of the
            Magnus expansion \cref{eq:Mag} is \small
            \begin{equation}
            \label{eq:gradfree}
            \tag{S3}
            \ee^{3 \ii h^{-2} \ve^{-1} \MM{q}^\top\MM{x}} \ee^{- 6  h^{-2} \MM{p}^\top \nabla} \ee^{T_2+(- \ii h \ve^{-1}  \widetilde{V} -6 \ii h^{-2} \ve^{-1} \MM{q}^\top\MM{x} + \tilde{c})}  \ee^{- 6  h^{-2} \MM{p}^\top \nabla} \ee^{3 \ii h^{-2} \ve^{-1} \MM{q}^\top\MM{x}} ,
            \end{equation}
             \normalsize where
            \begin{equation}
            \label{eq:ctilde}
            \tilde{c} = c - \ii \ve^{-1} \MM{q}^\top \MM{r}.
            \end{equation}
             In the context of \cref{eq:common},
            \begin{align}
            \nonumber T_3 &= T_2 = \ii h \ve \Delta   -\tilde{\MM{s}}^\top \nabla, & W_3 &= - \ii h \ve^{-1}  \widetilde{V} - L_3 + \tilde{c}, \\
            \label{eq:parts3} C_3 &= C_2 = - 12  h^{-2} \MM{p}^\top \nabla, & \ L_3 &= 6 \ii h^{-2} \ve^{-1} \MM{q}^\top\MM{x},
            \end{align}
             where $W_3$ is a slight perturbation of $W_2=W_1$ but is still
            a function (not a differential operator), while $C_3$ and $T_3$
            are identical to $C_2$ and $T_2$ defined in \cref{eq:parts2},
            respectively.

            \begin{remark}
            Once again, due to \cref{rmk:combineCT}, the exponential of $C_3$ can be
             combined with the exponential of $\ee^{T_3}$. Thus, the additional expense
             compared to a sixth-order scheme for time-independent potentials is only
             due to $\ee^{L_3}$, which is very inexpensive to compute. This is the
             marginal additional cost that we require in order to avoid computation of
             $\nabla V_0$ (in comparison to \cref{eq:commfree} described in
             \cref{sec:eliminatecomm}).
            \end{remark}

            \subsubsection{Derivation of the scheme}
            We start by attempting to express
            \begin{equation}\label{eq:derivationgradfree}
            \ee^{\Theta_4(t+h,t)} = \ee^{\frac12 L_3}\ee^{\frac12 C_2 } \ee^{T_2+W_3}\ee^{\frac12 C_2}\ee^{\frac12 L_3} + \OO{h^7}, \end{equation}
            where $\tilde{c},\gamma \in \BB{C}$ in $W_3$ and $L_3$,
            \[ L_3 = \gamma h^{-2} \MM{q}^\top\MM{x} = \OO{h^3},\]
            have to be determined.

            {\bf First application of sBCH.} We proceed by first expressing
            \[ \ee^{\frac12 C_2 } \ee^{T_2 + W_3}\ee^{\frac12 C_2 } = \exp(\mathrm{sBCH}(C_2,T_2 + W_3)), \]
             using the sBCH formula \cref{eq:sBCH} with $A= C_2$ and $B =
            T_2 + W_3$. Once again, expanding to grade three suffices due
            to the $\OO{h^3}$ size of $A=C_2$. In fact, the grade three
            commutator $[[B,A],A]$ can also be discarded for this very
            reason. For the only remaining commutator, using \cref{eq:ids}
            and the fact that $T_2$ commutes with $C_2$ the inner
            commutator $[B,A]$ reduces to
            \[[W_3,C_2] = -12 \ii h^{-1} \ve^{-1}  [\MM{p}^\top \nabla, \widetilde{V}] - 12 \gamma  h^{-4} [\MM{p}^\top \nabla,\MM{q}^\top\MM{x}] = -12 \ii h^{-1} \ve^{-1} \MM{p}^\top (\nabla V_0) + \alpha, \]
             where $\alpha$ is a scalar. Since this  commutator reduces to
            a function, it commutes with $W_3$ and, due to the $\OO{h^3}$
            size of $\MM{s}$, its commutator with $-  \tilde{\MM{s}}^\top
             \nabla$ is $\OO{h^7}$. Thus,
            \[ -\Frac{1}{12} [[B,A],B] = [\Delta, \MM{p}^\top (\nabla V_0)] + \OO{h^7} = C_1 + \OO{h^7}. \]
             We conclude
            \[\ee^{\frac12 C_2 } \ee^{T_2 + W_3}\ee^{\frac12 C_2 } = \exp(T_2 + W_3 + C_2 + C_1) + \OO{h^7}.\]

            {\bf Second application of sBCH.} In the second step, we
            express the right hand side of
            \[ \ee^{\frac12 L_3}\ee^{\frac12 C_2 } \ee^{T_2+W_3}\ee^{\frac12 C_2}\ee^{\frac12 L_3} = \ee^{\frac12 L_3}\exp(T_2 + W_3 + C_2 + C_1)\  \ee^{\frac12 L_3} + \OO{h^7},\]
             in terms of single exponential via sBCH. In this application
            we have $A = L_3$ and $B =T_2 + W_3 + C_2 + C_1$. Since $A =
            \OO{h^3}$ due to \cref{eq:powersofh}, $C_2 = \OO{h^3}$ and $C_1
            = \OO{h^5}$, any grade three commutator where more that one of
            these appear can be discarded, up to $\OO{h^7}$ accuracy. Once
            again, the only relevant commutator in the sBCH is $[[B,A],B]$,
            the relevant part of which is
            \[ [[B,A],B] = [[L_3,T_2+W_3],T_2+W_3] + \OO{h^7} = [[L_3,\ii h \ve \Delta],- \ii h \ve^{-1}  \widetilde{V}] + \OO{h^7},\]
             where we have used the fact that $L_3$ and $W_3$ commute,
            terms involving $\tilde{\MM{s}}^\top \nabla$ are $\OO{h^7}$ and
            that $[[L_3,\Delta],\Delta]$ vanishes due to \cref{eq:ids}.
            Thus,  this commutator simplifies to
            \[ [[B,A],B] = - \gamma [[\Delta, \MM{q}^\top\MM{x}],\widetilde{V}] =  - 2 \gamma [ \MM{q}^\top \nabla,\widetilde{V}] = - 2 \gamma \MM{q}^\top (\nabla V_0 + \MM{r}),\]
             using \cref{eq:ids}. Overall, in this step we find
            \[ \mathrm{sBCH}(A,B) = T_2 + W_3 + C_2  + C_1 + L_3+ \Frac16 \gamma \MM{q}^\top (\nabla V_0 + \MM{r}) + \OO{h^7}. \]
             This is identical to $\Theta_4 = T_1+W_1+C_1$ under the choice
            $\gamma = 6\ii \ve^{-1}$ and $\tilde{c} = c - \ii \ve^{-1}
            \MM{q}^\top \MM{r}$ since
            \[T_2+ C_2 = T_1, \qquad W_3 + L_3 + \ii \ve^{-1} \MM{q}^\top \nabla V_0 + \ii \ve^{-1} \MM{q}^\top \MM{r} = W_1,\]
             from \cref{eq:parts1,eq:parts2,eq:parts3}. Thus,
            \cref{eq:gradfree} approximates the exponential of the Magnus
            expansion \cref{eq:Mag} up to sixth-order accuracy.

\subsection{Approximation of the inner exponential}
\label{sec:inner} A wide range high-accuracy methods have been developed
            over the years
             \cite{mclachlan02sm,blanes06sso,blanes08sac,Omelyan}. Most of
            these can be readily employed for the approximation of
            $\exp(T+W)$ that is required in the schemes
            \cref{eq:stranglanczos,eq:commfree,eq:gradfree}, which share
            the common structure \cref{eq:common}.

            As we will discuss more concretely in
            \cref{sec:implementation}, the computation of the remaining
            exponentials in \cref{eq:common} can be very inexpensive under
            appropriate spatial discretisation strategies. In this sense,
            \cref{eq:common} is an effective strategy for extending
            existing methods for approximating $\exp(T+W)$ to the case of
            laser potentials at a very low cost.

            The choice of a method for approximating the inner exponential
            could be governed by a need for smaller error constants, better
            performance for large time steps and fewer exponentials, among
            other concerns. To describe concrete examples, in the following
            subsections we consider two different categories of splitting
            methods: a (i) classical splitting which only involves
            exponentials of $T$ and $W$ and (ii) a compact splitting where
            the number of exponential stages is reduced by utilising the
            gradient of $W$.

            The advantages of the proposed extension are not limited to
            splitting methods, however, and are equally applicable to
            situations where other approaches such as Lanczos
            approximation  \cite{parklight86} and Chebychev
            approximation  \cite{talezer84afs} for approximating $\exp(T+W)$
            prove more effective. As a concrete example we consider (iii)
            Lanczos approximation of the exponential since it is also
            relevant to other parts of the text.

            \subsubsection{Classical splittings}
            \label{sec:classical} As the first example of a splitting for
            $\ee^{T+W}$, we use the 15-stage sixth-order splitting
            specified by eqs. (84) and (85) in section~4.1.2 of Omelyan
            {\em et al.}  \cite{Omelyan},
            \begin{equation}\label{eq:classical}
            \tag{OMF85}
              \ee^{a_1 T} \ee^{b_1 W}  \ee^{a_2 T} \ee^{b_2 W} \ee^{a_3 T} \ee^{b_3 W}  \ee^{a_4 T}  \ee^{b_4 W} \ee^{a_4 T}  \ee^{b_3 W}  \ee^{a_3 T} \ee^{b_2 W} \ee^{a_2 T} \ee^{b_1 W}  \ee^{a_1 T},
            \end{equation}
             where
            \begin{align*}
            a_1 &= -1.0130879789171747, \quad & b_1 &= 0.00016600692650009894,\\
            a_2 &= 1.1874295737325427, \quad & b_2 &= -0.3796242142637736,\\
            a_3 &= -0.018335852096460590, \quad & b_3 &= 0.6891374118518106,\\
            a_4 &= 0.3439942572810926, \quad & b_4 &= 0.3806415909709257.
            \end{align*}
             Combining this with the outer exponents in any of the
            splittings \cref{eq:stranglanczos,eq:commfree,eq:gradfree},
            fully describes a concrete example of a sixth-order scheme for
            time-dependent potentials.

            In this way, any sixth-order classical splitting can be
            combined with any of the three approaches presented in this
            paper in a very straightforward way. Another alternative to
            \cref{eq:classical}, for instance, is described by eqs. (82)
            and (83) in section~4.1.1 of Omelyan {\em et
            al.}  \cite{Omelyan}, where the outermost exponentials are of $W$
            instead of $T$. This will be denoted by OMF83.

            \begin{remark}Note that in practice one might find that a fourth-order splitting with
            low error constant performs just as well for the approximation
            of $\exp(T+W)$, especially for large time steps. For this
            purpose, we will also consider the use of the fourth-order
            schemes OMF71 and OMF80, described by eqs. (63) and (71), and
            eqs. (72) and (80), respectively, in Omelyan {\em et
            al.}  \cite{Omelyan}.
            \end{remark}

            \subsubsection{Compact splittings}
            \label{sec:compact} Another concrete example results from using
            the 11-stage compact splitting given by eqs. (72) and (76) in
            section~3.5.2 of Omelyan {\em et al.}  \cite{Omelyan},
            \begin{equation}
            \label{eq:compact}
            \tag{OMF76}
            \ee^{a_1 T}\ee^{b_1 W + c_1 U} \ee^{a_2 T} \ee^{b_2 W + c_2 U}
            \ee^{a_3 T} \ee^{b_3 W + c_3 U}
            \ee^{a_3 T} \ee^{b_2 W + c_2 U} \ee^{a_2 T} \ee^{b_1 W + c_1 U} \ee^{a_1 T},
            \end{equation}
             where
            \begin{align*}
            a_1 &= 0.1097059723948682,   &b_1 &= 0.2693315848935301,    &c_1 &= 0.0008642161339706166,\\
            a_2 &= 0.4140632267310831,   &b_2 & = 1.1319803486515564,   &c_2 &= -0.01324638643416052,\\
            a_3 &=\Frac12-(a_1+a_2),     &b_3 & = 1-2(b_1+b_2),         &c_3 &= 0,
            \end{align*}
             and
            \[ U = -[[T,W],W],\]
             is a commutator of $T$ and $W$. An alternative with leading W,
            OMF65, is described by eqs. (63) and (65) in section~3.5.1 of
            Omelyan {\em et al.}  \cite{Omelyan}.

            In the case of the first splitting \cref{eq:stranglanczos}, $T$
            and $W$ are given by \cref{eq:parts1},
            \[T_1 =  \ii h \ve \Delta   - \MM{s}^\top \nabla, \quad W_1 =  - \ii h
            \ve^{-1} \widetilde{V} + \ii \ve^{-1} \MM{q}^\top (\nabla V_0)
             +c.\] Using \cref{eq:ids},
            \[ [[\ii h \ve \Delta, W_1],W_1] = 2 \ii h \ve (\nabla W_1)^2, \quad [[- \MM{s}^\top \nabla, W_1],W_1] = -[\MM{s}^\top (\nabla W_1),W_1] =0.\]
             so that $U = - \ii h \ve (\nabla W_1)^2$. This term possesses
            the same structure as $W_1$ (i.e. it is a function, not a
            differential operator) and combining it with $W_1$ in the
            splitting is a sensible approach.

            Additional care is required here, however, since the
            computation of $\nabla W_1$ can be very problematic  due to the
            presence of the $\MM{q}^\top (\nabla V_0)$ term in $W_1$. This
            would normally result in a need of the mixed derivatives,
            $\nabla(\MM{q}^\top (\nabla V_0))$, making computation very
            expensive. However, due to \cref{eq:powersofh}, $\MM{q}$ scales
            as $\OO{h^5}$ and, since $W_1 = \OO{h}$, the $\MM{q}^\top
            (\nabla V_0)$ term makes an $\OO{h^7}$ contribution to $U = -
            \ii h \ve (\nabla W_1)^2$ and can be ignored. Effectively, it
            suffices to use
            \begin{equation}\label{eq:Ucompact}U =  2\ii h^3 \ve^{-1} (\nabla
            \widetilde{V})^2.\end{equation}

            $U$ turns out to be the same in the case of \cref{eq:commfree} since $W_2 = W_1$, and $T_2$ differs from $T_1$ only in the $\MM{s}^\top \nabla$ term (see \cref{eq:parts2}), which does
            not contribute to $U$ (up to order six) due to the size.

%
%
%

            \subsubsection{Choice of classical vs compact splittings}
            In the case of fourth-order methods of Iserles {\em et
            al.}  \cite{IKS18cpc}, $\nabla \widetilde{V} = \nabla V_0 +
            \MM{r}(t,h)$ needs to be computed only in the case of compact
            splittings, and is not required for the classical splitting
            \cref{eq:MaStBM}.

            In contrast, in the sixth-order case $\nabla V_0$  appears directly in the
            Magnus expansion and is required in the schemes
            \cref{eq:stranglanczos,eq:commfree}, even when using classical splittings for
            the inner exponential. In these cases, the use of $\nabla V_0$ in the compact
            splitting \cref{eq:compact} is not an additional expense. Thus compact
            splittings should be favoured for \cref{eq:stranglanczos,eq:commfree}.

            Note that the typically expensive part, $\nabla V_0$, needs to be computed
            only once, while the computation of $\MM{r}$ is inexpensive.

            Although, in principle, it is possible to utilise a compact splitting for the
            central exponent of \cref{eq:gradfree}, it involves re-introducing the
            gradient of the potential, $\nabla V_0$, which defeats the point of
            \cref{eq:gradfree}.

\subsubsection{Lanczos approximation of the exponential}
            \label{sec:lanczos} Splitting methods are by no means the only
            approximation strategy for the exponential of a
            time-independent Hamiltonian. Lanczos approximation for
            $\exp(\AAA)\MM{v}$, where $\AAA \in \GG{su}(N)$ is an $N \times
            N$ skew-Hermitian matrix, involves approximating
            $\exp(\AAA)\MM{v}$ in the $m$th Krylov subspace,
            \begin{displaymath}
              \MM{K}_m(\mathcal{A},\MM{v})=\CC{span}\,\{\MM{v},\mathcal{A}\MM{v},\mathcal{A}^2\MM{v},\ldots,\mathcal{A}^{m-1}\MM{v}\},\qquad m\in\BB{N}.
            \end{displaymath}
            Effectively, the exponential is approximated in an $m$-dimensional subspace, as
            \begin{equation}
              \label{eq:KrylovExp}
              \ee^{\mathcal{A}}\MM{v}\approx \mathcal{V}_m\ee^{\mathcal{H}_m}\mathcal{V}_m^*\MM{v},
            \end{equation}
            where $\mathcal{V}_m$ is an orthogonal basis of the Krylov
            subspace $\MM{K}_m(\mathcal{A},\MM{v})$ and $\mathcal{H}_m$ is
            a tridiagonal matrix, both of which are found through the
            Lanczos iteration process. The Lanczos iterations involve
            computing matrix-vector products of the form $\AAA \MM{v}$,
            which are required for generating the Krylov subspace, combined
            with an orthogonalization procedure.

            When $m \ll N$, the exponential of the $N \times N$ matrix
            $\AAA$ is approximated by the exponential of the $m \times m$
            matrix $\mathcal{H}_m$ and the approximation is very
            inexpensive. The cost is largely dominated by the $m$
            matrix-vector products of the form $\AAA \MM{v}$. This approach
            has been used effectively in theoretical chemistry for
            long  \cite{parklight86} and its various aspects have been well
            studied  \cite{golub96mc}.

            The cost of computing $\AAA \MM{v}$ naturally depends on the
            structure of $\AAA$. In the case of general Magnus
            expansions  \cite{mko99fla,blanes00iho}, the matrix $\AAA$
            involves nested commutators, leading to a large cost of $\AAA
            \MM{v}$. In the case of the simplified Magnus expansion
            \cref{eq:Mag}, this cost is lower since it features a single
            commutator. This cost is lower still when we utilise Lanczos
            approximation
             of $\exp(T+W)$ inside \cref{eq:common} since $\AAA = T+W$ is
            free of commutators. Thus, this {\em extension} of Lanczos
            approximation may prove less expensive than other
            Magnus--Lanczos methods.

            The number of Lanczos iterations, $m$, which dictate the cost,
            naturally dictate the accuracy of the approximation
            \cref{eq:KrylovExp} as well. A tight error bound for the case
            of skew-Hermitian matrix is available  \cite{hochbruck97ksa},
            \begin{equation}
              \label{eq:LubichBound}
              \norm{2}{\ee^{\mathcal{A}}\MM{v}-\mathcal{V}_m\ee^{\mathcal{H}_m}\mathcal{V}_m^*\MM{v}} \leq 12\ee^{-\norm{}{\AAA}^2/(4m)} \left(\frac{\ee\norm{}{\AAA}}{2m}\right)^m,\qquad m\geq\norm{}{\AAA}.
            \end{equation}
            \begin{figure}[h]
                \centering
                \includegraphics{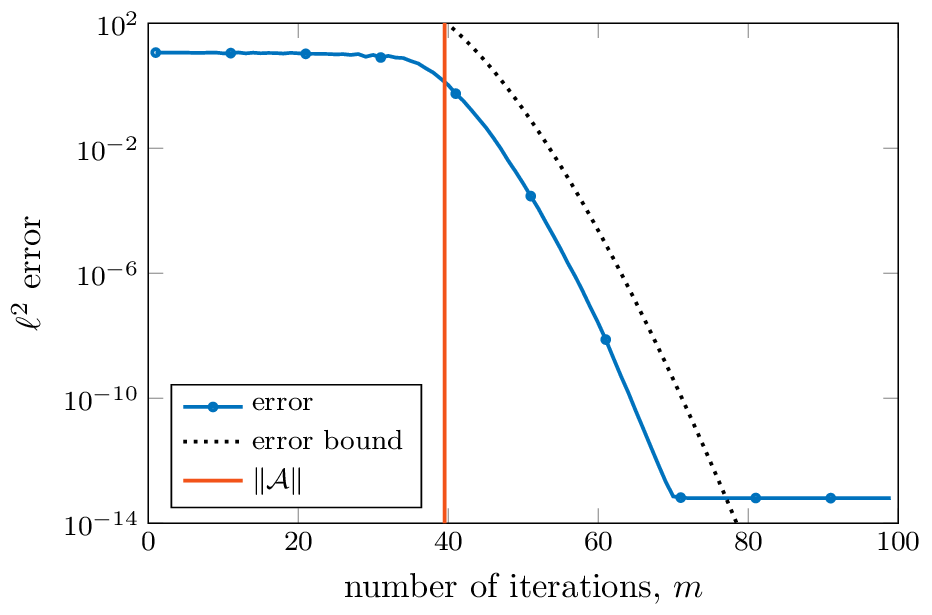}
                \caption{The Lanczos approximation to the matrix exponential for a random skew-Hermitian matrix $\AAA$ starts converging rapidly after $m \geq \norm{}{\AAA}$.}
                \label{fig:krylov}
            \end{figure}

            In particular, it has been well documented that the side
            condition $m\geq\norm{}{\AAA}$ is not an artificial imposition.
            In practice, Lanczos approximation of the exponential does not
            display superlinear convergence till the number of iterations
            have exceeded the spectral radius of the matrix $\AAA$ (see
            \cref{fig:krylov}).

\section{Implementation}
\label{sec:implementation}
\subsection{Approximation of integrals}
\label{sec:integrals}
            Depending on $\MM{e}$, analytic expressions for the integrals
            $\MM{r},\MM{s},\MM{q},\MM{p}$ and $c$ appearing in \cref{eq:Mag}, might be
            available. In the absence of analytic expressions, various quadrature methods
            can be utilised. For instance, if the quadrature weights and knots over
            $[0,h]$ are given by $w_1,\ldots,w_k$ and $\zeta_1,\ldots,\zeta_k$,
            respectively, we can approximate

            \begin{equation} \label{eq:quadrature} \Int{\zeta}{0}{h}{f(h,\zeta) \MM{e}(t+\zeta)} \approx
            \sum_{j=1}^k w_j f(h,\zeta_j) \MM{e}(t+\zeta_j),\end{equation}
             as usual. The nested integral in $c$, \cref{eq:c}, can be approximated as
             \begin{equation} \label{eq:quadraturenested}
            \Int{\zeta}{0}{h}{\zeta
            \MM{e}(t+\zeta)^\top \Int{\xi}{0}{\zeta}{\MM{e}(t+\xi)}}
            \approx \sum_{i=1}^k \sum_{j=1}^k \widetilde{w}_{ij} \zeta_i \MM{e}(t+\zeta_i)^\top \MM{e}(t+\zeta_j),\end{equation}
            where the weights,
            \[ \widetilde{w}_{ij} = \Int{\zeta}{0}{h}{
            \ell_i(t+\zeta) \Int{\xi}{0}{\zeta}{\ell_j(t+\xi)}},\]
            are found by substituting the Lagrange interpolating
            polynomial for $\MM{e}$,
            \[ \sum_{j=1}^k \ell_j(t+\zeta) \MM{e}(t+\zeta_j), \quad \ell_j(t+\zeta_i) = \delta_{ij}, \]
            in the integral. For order six accuracy, for instance, Gauss--Legendre quadrature with three knots suffice
            for non-oscillatory potentials. In this case, \footnotesize
            \[ \MM{w} = \Frac{h}{18}(5,8,5),\quad \MM{\zeta} = \Frac{h}{2}\left(1-\sqrt{3/5},\ 1,\ 1+\sqrt{3/5}\right),\]
            \normalsize and \footnotesize
            \[ \widetilde{\MM{w}} =  \frac{h^2}{648}
            \left(
              \begin{array}{ccc}
                25 & 40 - 12 \sqrt{15} & 25 - 6 \sqrt{15} \\
                40 + 12 \sqrt{15} & 64 & 40 - 12 \sqrt{15} \\
                25 + 6 \sqrt{15} & 40 + 12 \sqrt{15} & 25 \\
              \end{array}
            \right).\]  \normalsize

\subsection{Computation of exponentials}
\label{sec:expXYZ}
            The evaluation of exponentials of the modified kinetic and potential terms,
            $T_k$ and $W_k$, should be no more costly than the exponentiation of the
            Laplacian and the potential that are routinely employed in a sixth-order
            splitting scheme for time-independent potentials.

            In \cref{eq:commfree} and \cref{eq:gradfree}, $C_3=C_2$ shares the structure of $T_k$s, and in \cref{eq:gradfree} $L_3$
            has the same structure as $W_k$s. Consequently, the computation of their exponentials is not exceptionally problematic either.
            The commutator term in \cref{eq:stranglanczos}, $C_1$, however, requires a different strategy and
            is exponentiated via Lanczos iterations.

            \subsubsection{Exponentiating the modified potential terms -- $W_k$s and $L_3$}
            \label{sec:Y} Under spectral collocation on an equispaced grid,
            the term $W_1$ in \cref{eq:stranglanczos} discretises to a
            diagonal matrix,
            \[ W_1 \leadsto \DDD_{W_1} = \DDD_{- \ii h \widetilde{V} + \ii \ve^{-1} \MM{q}^\top (\nabla V_0) +c}\ , \]
             where $\leadsto$ denotes discretisation and $\DDD_{f}$ is a
            diagonal matrix with the values of $f$ on the grid points along
            its diagonal. The exponential of this is evaluated directly in
            a pointwise fashion,
            \[ \exp(\DDD_f) = \DDD_{\exp(f)}. \]
             The same holds true for the $W_k$s in \cref{eq:commfree} and
            \cref{eq:gradfree}, $L_3$ in \cref{eq:gradfree} and for
            exponentials of the form $\exp(b W + c U)$ in the splitting
            \cref{eq:compact}.

            \subsubsection{Exponentiating the modified kinetic terms -- $T_k$s and $C_2$}
            \label{sec:X} For the modified kinetic terms, we note that
            differentiation matrices are circulant and are diagonalised via
            Fourier transforms,
            \[ \partial_{x}^k \leadsto \MM{D}_{k,x} = \FFF_{x}^{-1} \DDD_{c_{k,x}} \FFF_{x}, \]
             where $\FFF_x$ is the Fourier transform in $x$ direction,
            $\FFF^{-1}_x$ is the inverse Fourier transform,
            $\DDD_{c_{k,x}}$ is a diagonal matrix and the values along its
            diagonal, $c_{k,x}$, comprise the symbol of the $k$th
            differentiation matrix, $\MM{D}_{k,x}$.

            In two dimensions, for instance, the exponential of the
            Laplacian term, $\ii h \ve \Delta$, is routinely computed in
            exponential splitting schemes for \schr equation with
            time-independent potentials as
            \[\ee^{\ii  h \ve \Delta} v \leadsto \FFF_{x}^{-1} \DDD_{\exp(\ii h \ve c_{2,x})} \FFF_{x} \FFF_{y}^{-1} \DDD_{\exp(\ii  h \ve c_{2,y})} \FFF_{y} \MM{v}, \]
             using four Fast Fourier Transforms (FFTs).

            Note that $\FFF_{y}\MM{v}$ is implemented in {\sc Matlab} as
            \texttt{fft(v,[],2)} for $\MM{v}$ discretised over an \texttt{ndgrid} domain,
            and $\DDD_{\MM{a}} \MM{v}$ is simply implemented as \texttt{a.*v}. Here
            $c_{2,x}$ and $c_{2,y}$ are the symbols of the differentiation matrices
            corresponding to $\dx^2$ and $\dy^2$, respectively. For instance, using a
            Fourier spectral method on an $[-L,L]^2$ box, where $\partial_y^2 \ee^{\ii
            \pi (j x + k y)/L} = (- \pi^2 k^2 /L^2) \ee^{\ii \pi (j x + k y)/L}$, we
            choose $(c_{2,y})_{j,k} = - \pi^2 k^2/L^2$.

            Using the same technique, we can compute the exponential of
            $T_1$ in \cref{eq:stranglanczos},
            \begin{align*}\ee^{\ii  h \ve \Delta -\frac12 \MM{s}(t,h)^\top
                \nabla} v &=
            \ee^{\ii h \ve \partial_{x}^2 - \frac12 s_x \partial_{x}}  \ee^{\ii h \ve
            \partial_{y}^2 - \frac12 s_y \partial_{y}} v \\
            &\leadsto \FFF_x^{-1} \DDD_{\exp( \ii h \ve c_{2,x} - \frac12
            s_x c_{1,x})}
             \FFF_x \FFF_y^{-1} \DDD_{\exp( \ii h \ve  c_{2,x} - \frac12
            s_y c_{1,x})}
             \FFF_y \MM{v},
            \end{align*}
             where $\MM{s} = (s_x,s_y)$, without any additional FFTs. The
            same observations apply to $T_k$s and $C_k$s in
            \cref{eq:commfree} and \cref{eq:gradfree}.

            \subsubsection{Exponentiating the commutator term -- $C_1$} \label{sec:Z}
             Unlike $T_k$s and $W_k$s,
             \[C_1 =\left[\Delta, \MM{p}^\top
            (\nabla V_0)\right],\]  which appears in
            \cref{eq:stranglanczos}, does not possess a structure that
            allows for direct exponentiation. However, the spectral radius
            of $C_1$ upon discretisation,
            \[C_1
            \leadsto \MM{C}_1 = \left[\sum_{j=1}^n \MM{D}_{2,x_j}\ ,\
            \DDD_{\MM{p}^\top (\nabla V_0)}\right],\]  is very small since
            $\MM{p} = \OO{h^5}$. Since $\MM{D}_{k,x}$ scales as
            $\OO{(\Delta x)^{-k}}$,
            \[ \rho(\MM{C}_1) = \norm{2}{\MM{C}_1} = \OO{h^5 (\Delta x)^{-2}} \norm{}{\nabla V_0}, \]
             assuming that we use the same spatial resolution in all
            directions.

            We can improve upon the estimate of the spectral radius
            further  \cite{IKS18jcp,bader14eaf}, whereby we find
            $\rho(\MM{C}_1) = \OO{h^5 (\Delta x)^{-1}} \max_{i,j \in
            \{1,\ldots,n \} } \norm{}{\partial_{x_i} \partial_{x_j} V_0}$.
            This improvement is notable in the semiclassical regime where a
            spatial resolution of $\Delta x = \OO{\ve}$ is necessitated by
            the highly oscillatory solution.

            This observation makes Lanczos iterations a very appealing
            candidate for the exponentiation of $\MM{C}_1$. As noted in
            \cref{sec:lanczos}, these methods feature a superlinear
            accuracy
             once the number of Lanczos iterations has exceeded the
            spectral radius of the exponent. In practice, in the case of
            $\MM{C}_1$, we find ourselves in the regime of superlinear
            accuracy of Lanczos iterations almost immediately and even a
            single Lanczos iteration seems to be giving us very good
            results.

            Each Lanczos iteration involves the computation of matrix-vector product of
            the form $\MM{C}_1 \MM{v}$, which can be computed as
            \begin{eqnarray*}
            \MM{C}_1 \MM{v} &=&  \sum_{j=1}^n \FFF_{x_j}^{-1} \DDD_{c_{2,x_j}} \FFF_{x_j}
            \DDD_{\MM{p}^\top (\nabla V_0)} \MM{v} -  \DDD_{\MM{p}^\top (\nabla V_0)}
            \sum_{j=1}^n \FFF_{x_j}^{-1} \DDD_{c_{2,x_j}} \FFF_{x_j} \MM{v}
            \end{eqnarray*}
             via $4n$ FFTs. Since directions are independent, these can be
            parallelized. Alternatively, one may use four n-dimensional FFTs (such as
            {\sc Matlab}'s \texttt{fftn}).

\section{Numerical examples}
\label{sec:numerics} In this section we provide detailed numerical
        experiments for two one-dimensional numerical examples considered
        by Iserles {\em et al.}  \cite{IKS18cpc} for which accurate reference
        solutions have been obtained by brute force (using extremely fine
        time steps and spatial grids).The first of these examples is in the
        regime, $\ve_1=1$, and the second in the semiclassical regime of
        $\ve_2 = 10^{-2}$. In both cases, we impose periodic boundaries on
        the spatial domains and resort to spectral collocation for
        discretisation.

        In principle, the procedure extends in a straightforward way to
        higher dimensions via tensorisation of the periodic grids and we
        demonstrate the applicability of the approach using two and
        three-dimensional examples under $\ve_3=\ve_4=10^{-2}$. Lastly, we
        consider a Coloumb potential example from Schaeffer {\em et
        al.}  \cite{Schaefer} under $\ve_5=1$.

\begin{figure}[tbh]
            \includegraphics[width=13cm]{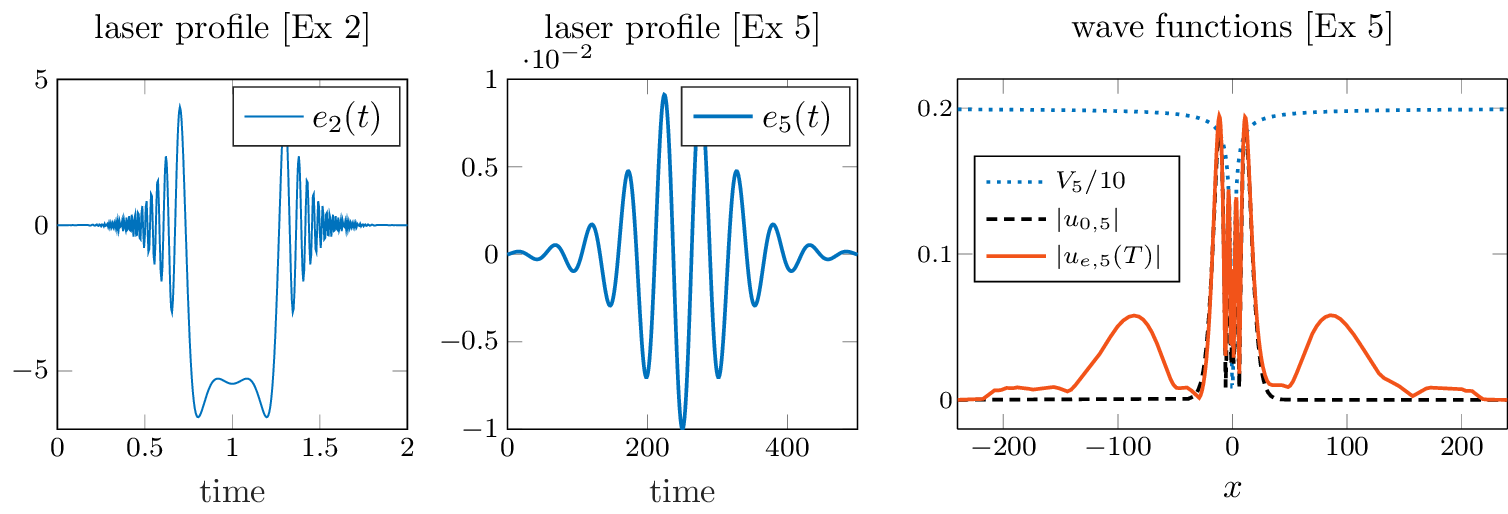}
        \caption{[Ex 2, 5] Laser profiles $e_2$ (left) and $e_5$ (centre). The fifth eigenfunction of the soft Coloumb potential $V_5$
        evolves to $u_{e,5}(T)$ at $t=T$ under the influence of $e_5$ (right). }
        \label{fig:uColoumb}
    \end{figure}

        In the first two examples, the initial conditions $u_{0,1}$ and
        $u_{0,2}$ are Gaussian wavepackets,
        \[u_{0,k}(x) = (\delta_k  \pi)^{-d/4}  \exp\left( (-(x - x_0)^2 )/ (2 \delta_k)\right),\quad x_0 = -2.5,\ k=1,2,\]
         with $d=1$ (for the one dimensional problems), $\delta_1 = 0.2$
        and $\delta_2 = 10^{-2}$ in the respective cases. These wavepackets
        are sitting in the left well of the double well potentials,

        \[V_{1}(x) = x^4 - 15 x^2 \quad \mathrm{and}\quad V_{2}(x) = \Frac15 x^4 - 2 x^2, \]
         respectively, which act as the choice of $V_0$ in the two
        examples. The time profile of the laser used here is
        \[ e_1(t) = \begin{cases} \sin(25 \pi t) & \quad t \in [\Frac35 n,\Frac35 n + \Frac{1}{25}],\quad n \geq 1,\\
                                  \sin(5 \pi t) & \quad t \in (\Frac35 n + \Frac{1}{25},\Frac35 n + \Frac{6}{25}],\quad n \geq 1,
                    \end{cases} \]
        and
        \[ e_2(t) = 10 \exp(-10 (t - 1)^2) \sin((500 (t - 1)^4 + 10)), \]
        respectively. The former is a sequence of {\em asymmetric sine
        lobes} while the latter is a highly oscillatory {\em chirped} pulse
        (\cref{fig:uColoumb} (left)). Such laser profiles are used
        routinely in laser control  \cite{AmstrupChirped}. Even more
        oscillatory electric fields often result from optimal control
        algorithms  \cite{MeyerOptimal,CoudertOptimal}.

        The spatial domain is $[-10,10]$ and $[-5,5]$ in the two examples,
        respectively, while the temporal domain is $[0,4]$ and $[0,\frac52]$,
        respectively.

\begin{figure}[tbh]
            \includegraphics[width=13cm]{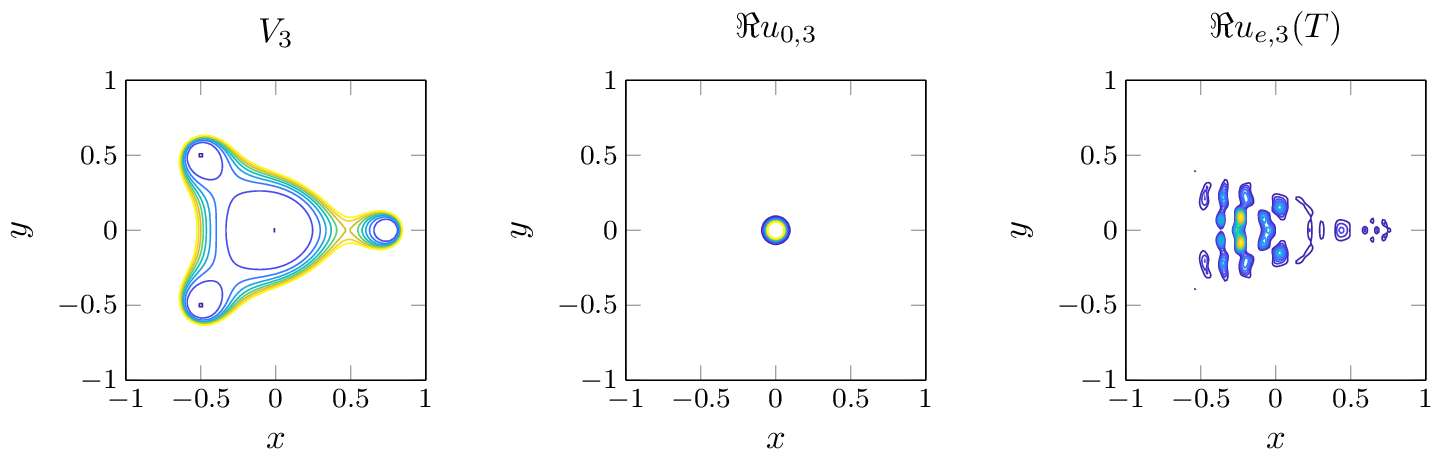}
        \caption{[Ex 3] The potential $V_{3}$ in two dimensions shown with contour lines (left); initial condition $u_{0,3}$  (centre) evolves to
        to $u_{e,3}(T)$ (right) under the influence of $e_3(t) = \Frac15 e_2(t) (1,0)^\top$. Contour lines for the  potential range from  $0$ to $2$, while those
        for (real part of the) wavefunctions range from  $0.25$ to $4$. The levels are separated by $0.25$ in all cases.
        Here, the semiclassical parameter is $\ve_3 = 1/100$. }
        \label{fig:uIF2D}
    \end{figure}

\begin{figure}[tbh]
            \includegraphics[width=13cm]{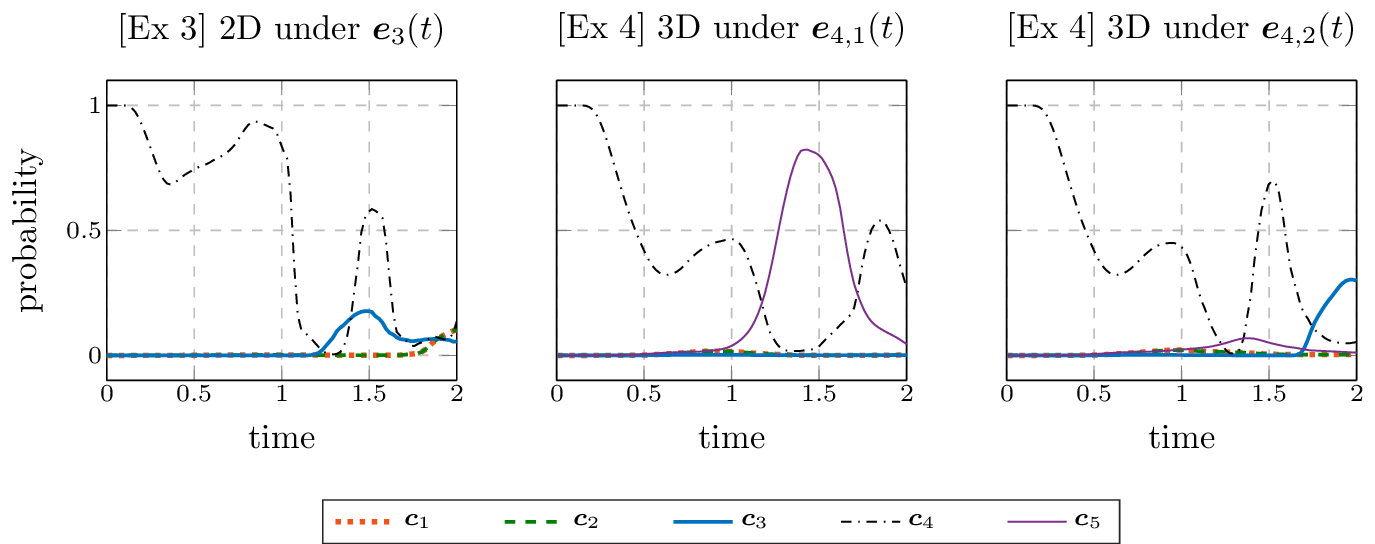}
        \caption{[Well Occupation] Probability of particle being in a radius $0.2$ from the centres $\MM{c}_j$ in two dimensions under the influence of $\MM{e}_3(t)$ (left) and in three dimensions under the influence of $\MM{e}_{4,1}(t)$ (centre) and $\MM{e}_{4,2}(t)$ (right). }
        \label{fig:targets}
        \end{figure}

        In the third and fourth examples, we consider the Gaussian
        wavepackets
        \[u_{0,k}(x) = (\delta_k  \pi)^{-d/4}  \exp\left( -\norm{}{\MM{x}}^2/ (2 \delta_k)\right), \quad k=3,4,\]
        with $d=2$ for a two-dimensional example (example 3) and $d=3$ for a three-dimensional example (example 4).
        We use $\delta_3 = \delta_4 = 10^{-3}$ and $\ve_3=\ve_4=10^{-2}$ in both cases.
        These wavepackets are sitting in the central wells of the potentials
        \[ V_3(\MM{x}) = 2500 \prod_{j=1}^4 \norm{}{\MM{x}-\MM{c}_j}^2, \qquad V_4(\MM{x}) = 40 \prod_{j=1}^5 \norm{}{\MM{x}-\MM{c}_j}^2,\]
        which are degree eight and degree ten polynomials in two and three dimensions, respectively, and where the wells are defined by the centres,
        \footnotesize
        \[  \left(
                     \MM{c}_1, \ldots, \MM{c}_4
                 \right) = \left(
                             \begin{array}{ccccc}
                               -0.5 & -0.5 & 1/\sqrt{2} & 0  \\
                               -0.5 & 0.5 & 0 & 0 \\
                             \end{array}
                           \right),
         \quad \left(
            \MM{c}_1, \ldots, \MM{c}_5
             \right) = \left(
                         \begin{array}{ccccc}
                           -0.5 & -0.5 & 0.75 & 0 & 0 \\
                           -0.5 & 0.5 & 0 & 0 & 0 \\
                           -0.5 & -0.5 & -0.5 & 0 & 0.75 \\
                         \end{array}
                       \right),
         \]
        \normalsize
        in the two cases, respectively. The spatial domain used is $[-1,1]^d$ and the temporal domain is $[0,2]$.
        In the two dimensional case (\cref{fig:uIF2D}) we consider the influence under a laser in the $x$ direction (\cref{fig:targets} (left)), while in the
        three dimensional case we consider lasers in two different directions of polarisation (\cref{fig:targets} (centre) and (right)),
        \[ \MM{e}_3(t) = \frac{e_2(t)}{5}(1,0)^\top, \qquad  \MM{e}_{4,1}(t) = \frac{e_2(t)}{5}(0,0,1)^\top, \quad \MM{e}_{4,2}(t) = \frac{e_2(t)}{5\sqrt{2}}(-1,0,1)^\top. \]

        Lastly, we consider the one-dimensional (soft) Coloumb potential
        example from Schaeffer {\em et al.}  \cite{Schaefer} as our fifth
        example (\cref{fig:uColoumb} (centre) and (right)), where $V_5(x) = 2\left(1-\Frac{1}{\sqrt{x^2+1}}\right)$.
        We take the (numerically determined) fifth eigenfunction of this
        potential as the initial condition, $u_{0,5}$, and consider its
        evolution under the influence of the laser profile
        \[ e_5(t) = -0.01\ \mathrm{sech}\left(\frac{t-250}{85}\right) \cos(0.12\,(t-250)),\]
        over a temporal domain $[0,500]$ and spatial domain $[-240,240]$.
        The temporal domain is halved and the potential scaled up by a
        factor of two to account for the fact that under our scaling
        $\Delta$ appears without a factor of $1/2$ in \cref{eq:schr}.

        The effective time-dependent potentials in our examples are
        \[ V_{e,k}(\MM{x},t) = V_{k}(\MM{x}) + \MM{e}_k(t)^\top \MM{x}, \quad k=1,2,3,4,5,\]
        with an additional index being used in the case of $k=4$ to
        distinguish the use of $\MM{e}_{4,1}$ or $\MM{e}_{4,2}$. Under the
        influence of these potentials, the initial conditions evolve to
        $u_{e,k}(T)$ at $t=T$. The probability of being within a radius of
        $0.2$ from the centres $\MM{c}_j$ in the two and three-dimensional
        examples under the influence of these lasers is shown in
        \cref{fig:targets}.

\begin{figure}[tbh]
            \centering
                \includegraphics{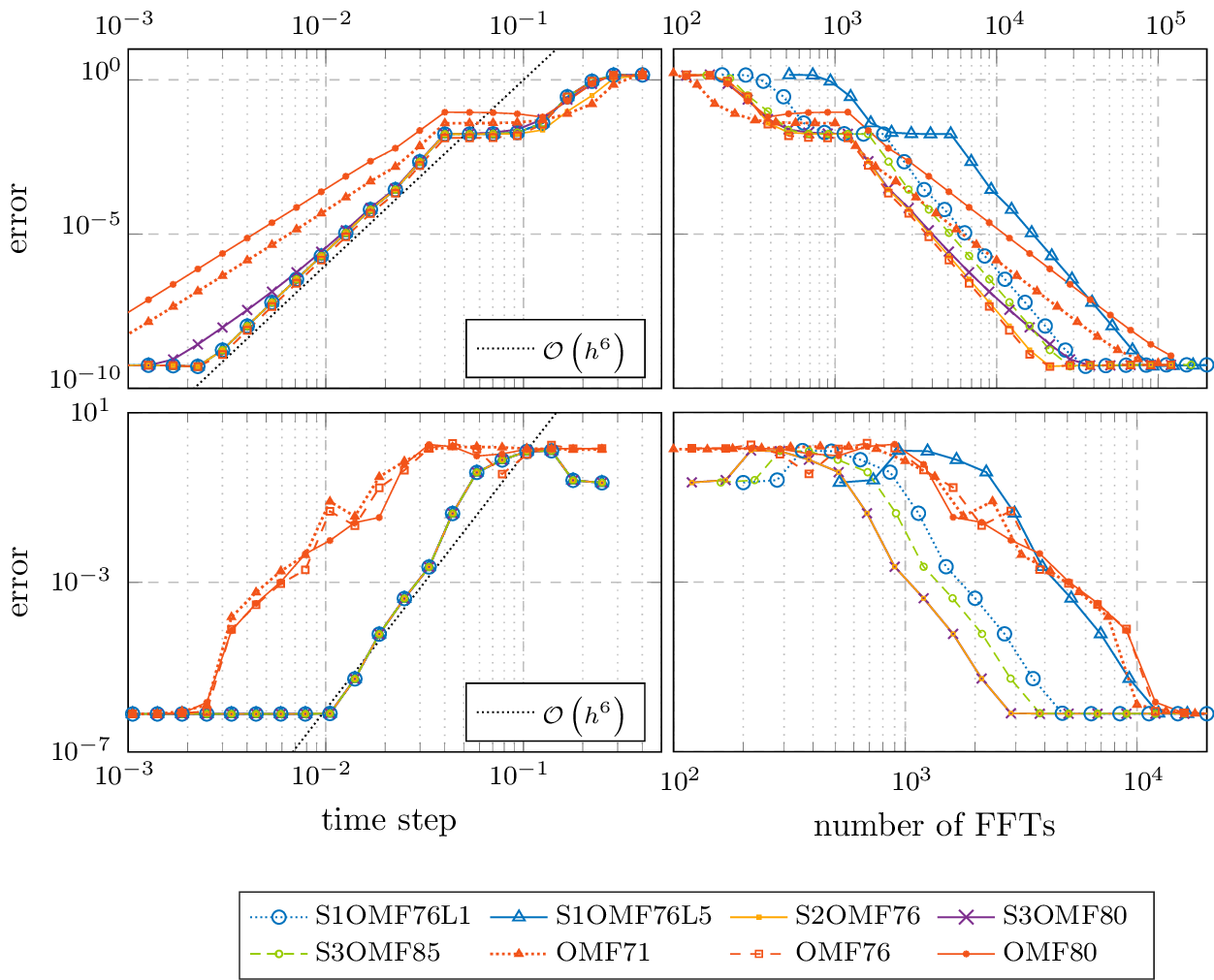}
                \caption{[Ex 1 (top), Ex 2 (bottom)] accuracy (left);
            efficiency (right).} \label{fig:errorSn}
        \end{figure}

\begin{figure}[tbh]
            \centering
                \includegraphics[width=13cm]{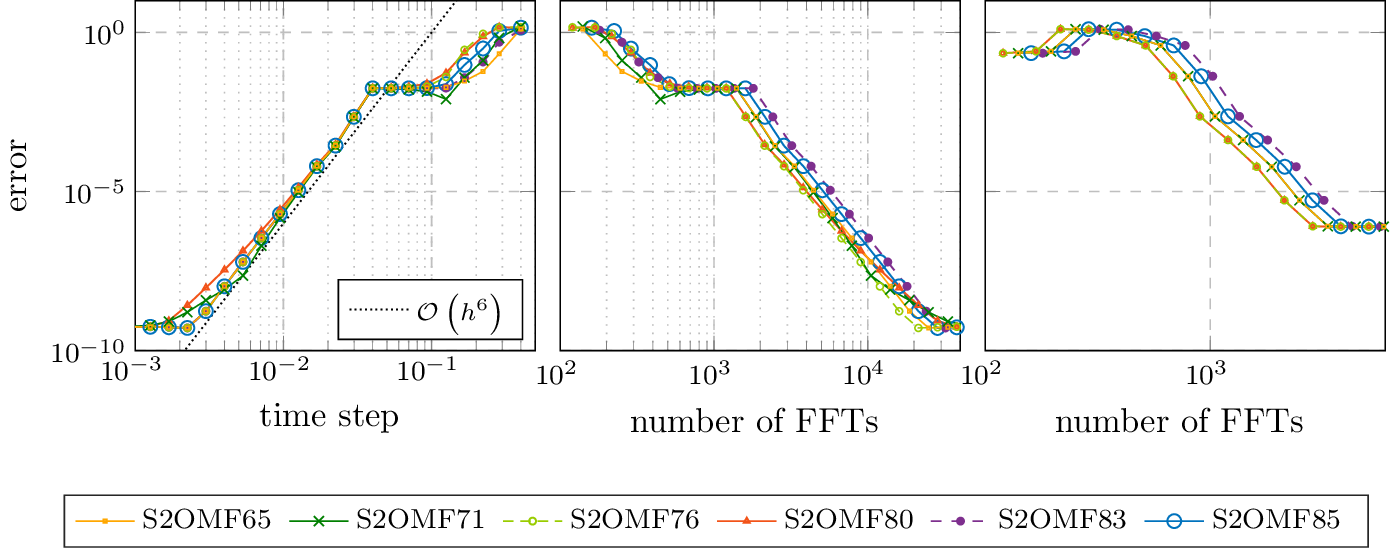}
                \caption{Accuracy (left) and efficiency for Ex 1 (centre); efficiency for Ex 2 (right).}
            \label{fig:errorS2OMFn}
        \end{figure}

         {\bf Methods.} In \cref{fig:errorS2OMFn,fig:errorSn} we display
        the comparisons of accuracy and efficiencies of several methods.
        The naming of methods is straightforward -- the combination of the
        $n$-th proposed scheme with the $k$-th splitting (the splitting
        described by the coefficients in eq ($k$)) of Omelyan {\em et
        al.}  \cite{Omelyan}  is labeled S$n$OMF$k$. For instance, the
        combination of \cref{eq:commfree} with \cref{eq:compact}, which
        describes a concrete scheme will be called S2OMF76, while the
        combination of \cref{eq:gradfree} with \cref{eq:classical} will be
        called S3OMF85. In the case of \cref{eq:stranglanczos}, we add the
        postfix L$m$ to denote the number of Lanczos iterations used for
        $C_1$.

        These methods are also compared against OMF$k$ (without a prefix of
        S$n$), which denote the time-ordered exponential splittings given
        by the coefficients in eq ($k$) of Omelyan {\em et
        al.}  \cite{Omelyan} that handle a time-dependent potential by
        advancing time along with the application of the Laplacian
          \cite{Goldstein2004}. These are an alternative to the proposed
        approach that were mentioned in the introduction.

        Commutator-free Lanczos based methods  \cite{alvermann2011hocm} are
        not compared against the proposed schemes since their relative
        ineffectiveness in the context of \cref{eq:schr} has already been
        studied  \cite{IKS18cpc}.

        In the first example $\ve_1=1$, we use $M_1=150$ spatial grid
        points while for second example, which features highly oscillatory
        behaviour in the solution due to the small semiclassical parameter
        $\ve_2=10^{-2}$, we use $M_2=1000$ spatial grid points (which,
        nevertheless, proves inadequate to achieve accuracies higher than
        $10^{-6}$). In the third and fourth we use $M=150$ points in each
        direction to keep computations manageable. This low spatial
        resolution limits the accuracy that can be achieved. In the final
        example, we use $M=768$ points for the Coloumb potential.

        {\bf Absorbing boundary.} We use an absorbing boundary of width
        $40$ using the procedure described in section~4.2 of Schaeffer {\em
        et al.}  \cite{Schaefer}. In particular, we need to compute the
        potential $V_{mod}$ which is a version of $V_{e,5}(x,t)$ that
        smoothly becomes flat in the absorbing boundary region (eq. (103)
        of Schaeffer {\em et al.}  \cite{Schaefer}). In principle this needs
        to be computed at each $t$. The problem is more easily resolved by
        computing the flattened versions $V_{5,mod}$ and $x_{mod}$ for
        $V_5$ and $x$ once, after which the flattened version of
        $V_{e,5}(x,t)$ can readily be computed as $V_{5,mod} + e_5(t)
        x_{mod}$ for any $t$. This directly leads to a flattened version of
        $\tilde{V}$ required in \cref{eq:gradfree}, to which we add the
        complex-valued absorbing potential from Schaeffer {\em et
        al.}  \cite{Schaefer} to complete the description.

        {\bf Reference solution.} In the first example the reference
        solution is obtained by using a sixth-order Magnus--Lanczos
        method \cite{BIKS}, while in the second example the reference
        solution is obtained using a sixth-order commutator-free Lanczos
        based method \cite{alvermann2011hocm}. In both cases we used $5000$
        spatial grid points and $10^6$ time steps. For the third and fourth
        examples, we generate a reference solution with $150$ grid points
        in each direction due to computational complexity, using Strang
        splitting with $2\times 10^6$ time steps for the third example and
        OMF80 with $10^4$ time steps for the fourth example. The reference
        solution for the fifth example is generated using OMF85 with
        $2\times 10^5$ time steps.

        {\bf Quadrature points.} In the first example, the integrals in our schemes
        are discretised via three Gauss--Legendre knots, while in the second example
        eleven Gauss--Legendre knots are used in order to adequately resolve the
        highly oscillatory potential. In contrast, the method OMF$k$ effectively use
        a fixed number of knots, dictated by the number of exponentials of the
        Laplacian. For instance, \cref{eq:compact} and OMF80 use five knots while
        OMF71 uses six knots.


        \begin{figure}[tbh]
            \includegraphics[width=13cm]{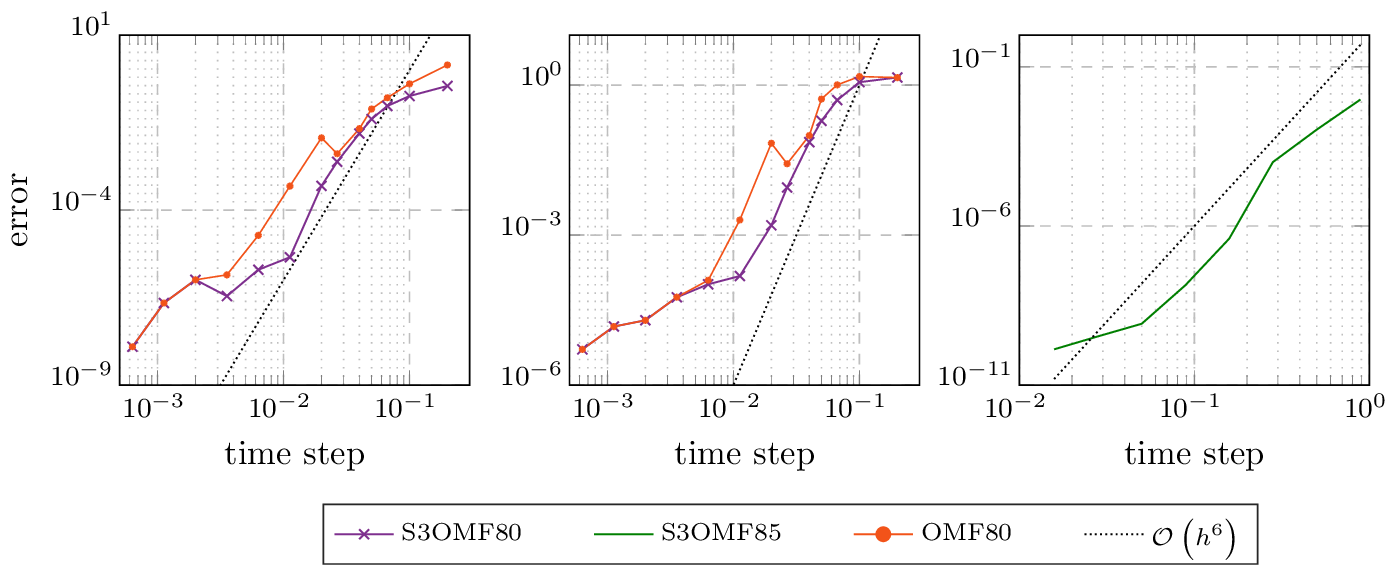}
        \caption{Accuracy for Ex 3 (left), Ex 4 under $\MM{e}_{4,1}(t)$ (centre) and Ex 5 (right).}
        \label{fig:errorND}
        \end{figure}

\section{Conclusions}
\label{sec:conclusions}
        We have presented three different strategies for
        easily extending existing sixth-order schemes for the \schr equation with
        time-independent potentials to the case of laser potentials under the dipole
        approximation. The overall schemes require, at most, one additional
        exponential, which leads to a very marginal increase in cost.

        Keeping the integrals $\MM{r},\MM{s},\MM{q},\MM{p}$ and $c$ intact in our
        schemes allows us flexibility in deciding a quadrature strategy at the very
        end. The advantage over the time-ordered exponential splittings OMF71, OMF76
        and OMF80, which sample the potential at fixed time knots, is evident from
        the numerical results for the second example (\cref{fig:errorSn} (bottom
        row)), where a highly oscillatory laser is involved.  Where a weaker field is
        involved (such as Examples 3 and 4 where $\norm{}{\MM{e}_3(t)} = \norm{}{\MM{e}_{4,1}(t)} = \norm{}{\MM{e}_{4,2}(t)} = \norm{}{\MM{e}_2(t)}/5$)
        and the solution is less oscillatory, the advantage
        may be less pronounced (\cref{fig:errorND} (left) and (centre)).

        Moreover, as mentioned in the introduction, preserving the integrals till the
        end allows us to use other quadrature methods which can be particularly
        helpful when the laser potential is only known at specific points (such as in
        control applications).

        The proposed methods are also effective when the potential is not
        highly oscillatory and can be sampled arbitrarily. As seen by
        numerical results for the first example (\cref{fig:errorSn} (top
        row)), the proposed schemes end up being very effective (nearly as
        accurate as the sixth-order methods) even when combined with OMF71
        and OMF80, which are fourth-order schemes with low error constants.
        This is in contrast to the direct use of OMF71 and OMF80 as
        time-ordered exponential splittings. A comparison of some fourth
        and sixth-order splittings from Omelyan {\em et al.} \cite{Omelyan}
        for the central exponent is presented in \cref{fig:errorS2OMFn}.

        Where the gradient of the potential, $\nabla
        V_0$, is available, we recommend using \cref{eq:commfree} in conjunction with
        compact splittings such as \cref{eq:compact}. Where $\nabla V_0$ needs to be
        avoided, we recommend using \cref{eq:gradfree} in conjunction with classical
        splittings such as \cref{eq:classical} or OMF80. The effectiveness of
        such a gradient free scheme, S3OMF85, is demonstrated for the soft Coloumb potential of Example 5 in
        \cref{fig:errorND} (right).

        We remind the reader that the splitting schemes of Omelyan {\em et al.} \cite{Omelyan} are
        among a myriad possible ways of approximating the innermost exponential, $\exp(T+W)$, in \cref{eq:common}.
        As mentioned in \cref{sec:inner}, for instance,
        alternative methods such as Lanczos approximation \cite{parklight86} and Chebychev
        approximation \cite{talezer84afs} can also be employed for this purpose.
        In contrast to a direct application of these techniques for exponentiating the Magnus expansion,
        the application to $\exp(T+W)$ would benefit from lower costs of matrix-vector products due to a commutator-free exponent.
        The primary focus of this manuscript, however, is to present a few general strategies (the splittings \cref{eq:stranglanczos,eq:commfree,eq:gradfree}) for extending existing schemes
        to the case of laser potentials, and the potential merits and limitations of following this approach in the context of various existing schemes have not been fully explored.


        While the numerical examples demonstrate the effectiveness of the
        schemes for laser potentials in one, two and three dimensions (in
        Cartesian coordinates), there are a range of issues and possible
        avenues of future research that our work raises:
        \begin{enumerate}
        \item {\em High dimensions.} For dimensions higher than three, tensorisation is not a viable strategy. Among
        the various approaches for truly high-dimensional problems, the extension of the proposed schemes to the case of
        Hagedorn wavepackets \cite{gradinaru14nm} is worth exploring since the additional terms in our schemes are at most
        quadratic in the momentum and position operators, $-\ii \nabla$ and $\MM{x}$.
        \item {\em Non-Cartesian coordinates.}
        While applicability of this approach to vibrational coordinates is unlikely, an approach for spherical coordinates
        using a mix of Fast Spherical Harmonic Transforms and Fast Fourier Transforms is being explored.
        \item {\em Matrix-valued potentials.} In the coherent control of a two-level atom \cite{shapiro03qc}, the potential becomes matrix-valued. The change in the algebraic nature of the problem
        necessitates the development of specialised splittings for this case.
        \end{enumerate}

\bibliographystyle{plain}
\bibliography{HighOrderLaser}

\appendix
\section{Scalar phase factors}
\label{app:scalars}
        \subsection{Simplification of $c_{3,1}$}
        \label{app:c31} For the simplification of the integrals in $c_{3,1}$, we
        define
        \begin{eqnarray*}
        I_1 & = & \Int{\zeta}{0}{h}{\zeta \MM{e}(t+\zeta)^\top
        \Int{\xi}{0}{\zeta}{\MM{e}(t+\xi)}},\\
        I_2 & = & \Int{\zeta}{0}{h}{ \left( \Int{\xi}{0}{\zeta}{\MM{e}(t+\xi)}
        \right)^2},\\
        I_3 & = & \Int{\zeta}{0}{h}{\MM{e}(t+\zeta)^\top
        \Int{\xi}{0}{\zeta}{\xi \MM{e}(t+\xi)}},\\
        \end{eqnarray*}
        and derive the following identities via integration by parts,
        \[ I_2 = - I_1 + I_3 + \left(\Int{\zeta}{0}{h}{\MM{e}(t+\zeta)}\right)^\top \left(\Int{\zeta}{0}{h}{(h-\zeta)\MM{e}(t+\zeta)}\right),\]
        and
        \[ I_3 = - I_1 +  \left(\Int{\zeta}{0}{h}{\MM{e}(t+\zeta)}\right)^\top \left(\Int{\zeta}{0}{h}{\zeta\MM{e}(t+\zeta)}\right),\]
        so that $I_2 = -2 I_1 +
         h\left(\Int{\zeta}{0}{h}{\MM{e}(t+\zeta)}\right)^2$, where
         $\MM{a}^2 = \MM{a}^\top \MM{a}$. Putting these together,
        \begin{align}
        \nonumber c_{3,1} &= \Frac{1}{6} \ii \ve^{-1}  \Int{\zeta}{0}{h}{\Int{\xi}{0}{\zeta}{
        (\MM{e}(t+\zeta)-\MM{e}(t+\xi))^\top \Int{\chi}{0}{\zeta}{ \MM{e}(t+\chi)
        }}}\\
        \label{eq:c31}  &=\Frac{1}{6} \ii \ve^{-1} (I_1 - I_2) =  \ii \ve^{-1} \left(\Frac{1}{2} I_1 - \Frac{1}{6}  h\left(\Int{\zeta}{0}{h}{\MM{e}(t+\zeta)}\right)^2\right).
        \end{align}

        \subsection{Simplification of $c_{3,2}$}
        \label{app:c32} For the simplification of $c_{3,2}$, we simplify
        \begin{eqnarray*}
        &&\Int{\zeta}{0}{h}{\MM{e}(t+\zeta)^\top
        \Int{\xi}{0}{\zeta}{\Int{\chi}{0}{\xi}{\MM{e}(t+\chi)}}} -
        \Int{\zeta}{0}{h}{\MM{e}(t+\zeta)^\top
        \Int{\xi}{0}{\zeta}{\Int{\chi}{0}{\xi}{\MM{e}(t+\xi) }}}\\
        & & \qquad = \Int{\zeta}{0}{h}{\MM{e}(t+\zeta)^\top
        \Int{\xi}{0}{\zeta}{(\zeta-\xi)\MM{e}(t+\xi) }} -
        \Int{\zeta}{0}{h}{\MM{e}(t+\zeta)^\top
        \Int{\xi}{0}{\zeta}{\xi\MM{e}(t+\xi) }}\\
        & & \qquad = I_1 - 2I_3 = 3I_1 - 2\left(\Int{\zeta}{0}{h}{\MM{e}(t+\zeta)}\right)^\top \left(\Int{\zeta}{0}{h}{\zeta\MM{e}(t+\zeta)}\right) .\\
        \end{eqnarray*}
        using \cref{eq:ibp1} under $n=0$. Consequently,
        \begin{align}
        \nonumber c_{3,2} &= \Frac12 \ii \ve^{-1} \Int{\zeta}{0}{h}{\MM{e}(t+\zeta)^\top
        \Int{\xi}{0}{\zeta}{\Int{\chi}{0}{\xi}{\MM{e}(t+\chi) - \MM{e}(t+\xi)}}} \\
        \label{eq:c32} &= \ii \ve^{-1} \left(\Frac32  I_1 -  \left(\Int{\zeta}{0}{h}{\MM{e}(t+\zeta)}\right)^\top \left(\Int{\zeta}{0}{h}{\zeta\MM{e}(t+\zeta)}\right) \right) .
        \end{align}
         Putting these together, $c(t,h)=c_{3,1}+c_{3,2}$,
        \begin{equation}
        \tag{\ref{eq:c}} c(t,h) =  \ii \ve^{-1} \left(2  I_1  - \left(\Int{\zeta}{0}{h}{\MM{e}(t+\zeta)}\right)^\top \left(\Int{\zeta}{0}{h}{\zeta\MM{e}(t+\zeta)}\right)   - \Frac{1}{6}   h \left(\Int{\zeta}{0}{h}{\MM{e}(t+\zeta)}\right)^2\right).
        \end{equation}

\end{document}